\documentclass[preprint,10pt]{elsarticle}

\usepackage{graphicx}%
\usepackage{multirow}%
\usepackage{amsmath,amssymb,amsfonts}%
\usepackage{amsthm}%

\usepackage{mathrsfs}%
\usepackage[title]{appendix}%
\usepackage{xcolor}%
\usepackage{textcomp}%
\usepackage{manyfoot}%
\usepackage{booktabs,mdframed}%

\usepackage{enumerate,courier}
\usepackage{subfigure,subfloat}
\usepackage{caption}
\usepackage{float}
\usepackage{graphicx}
\usepackage{epstopdf}
\usepackage{algorithm}
\usepackage{algpseudocode}

\usepackage {color}
\usepackage{epstopdf}
\usepackage{threeparttable}
\usepackage{setspace}
\usepackage{rotating}
\usepackage{varwidth}
\usepackage{caption}
\usepackage{hyperref} 
\usepackage{rotating}
\usepackage{caption}
\usepackage[margin=2.5cm]{geometry}	

\usepackage{enumitem}
\setlist[itemize]{itemsep=-1pt} 
\setlist[enumerate]{itemsep=-1pt} 

\hypersetup{
	colorlinks=true,
	citecolor=blue,
	linkcolor=red,
	urlcolor=green}
\newtheorem{theorem}{Theorem}
\newtheorem{proposition}[theorem]{Proposition}%

\numberwithin{theorem}{section}
\numberwithin{equation}{section}  
\newtheorem{assumption}{Assumption}[section]%
\newtheorem{remark}{Remark}[section]%
\newtheorem{lemma}{Lemma}[section]
\newtheorem{corollary}{Corollary}[section]
\newtheorem{definition}{Definition}[section]%

\begin{document}

\begin{frontmatter}



\title{Conjugate gradient methods without line search for multiobjective optimization} 


\author{Wang Chen\textsuperscript{a,b}, Yong Zhao\textsuperscript{c}, Liping Tang\textsuperscript{a}, Xinmin Yang\textsuperscript{a}}
\address[mymainaddress]{National Center for Applied Mathematics in Chongqing, Chongqing Normal University, Chongqing, 401331, China}
\address[mysecondaddress]{School of Mathematical Sciences, University of Electronic Science and Technology of China, Chengdu, 611731, China}
\address[mythirdaryaddress]{College of Mathematics and Statistics, Chongqing Jiaotong University, Chongqing, China}

\cortext[mycorrespondingauthor]{ \indent \indent  Email addresses: 
\href{mailto:chenwangff@163.com}{chenwangff@163.com} (Wang Chen),
\href{mailto:zhaoyongty@126.com}{zhaoyongty@126.com} (Yong Zhao),
\href{mailto:tanglipings@163.com}{tanglipings@163.com} (Liping Tang), \href{mailto:xmyang@cqnu.edu.cn}{xmyang@cqnu.edu.cn} (Xinmin Yang)}

\begin{abstract}
This paper addresses unconstrained multiobjective optimization problems where two or more continuously differentiable functions have to be minimized. We delve into the conjugate gradient methods proposed by Lucambio P\'{e}rez and Prudente (SIAM J Optim, 28(3): 2690--2720, 2018) for such problems. Instead of the Wolfe-type line search procedure used in their work, we employ a fixed stepsize formula (or no-line-search scheme), which can mitigate the pressure of choosing stepsize caused by multiple inequalities and avoid the computational cost associated with function evaluations in specific applications. 
The no-line-search scheme is utilized to derive the condition of Zoutendijk's type. Global convergence encompasses the vector extensions of Fletcher--Reeves, conjugate descent, Dai--Yuan, Polak--Ribi\`{e}re--Polyak and Hestenes--Stiefel parameters, subject to certain mild assumptions. Additionally, numerical experiments are conducted to demonstrate the practical performance of the proposed stepsize rule, and comparative analyses  are made with the multiobjective steepest descent methods using the Armijo line search and the multiobjective conjugate gradient methods using the Wolfe-type line search.
\end{abstract}

%


\end{frontmatter}



\section{Introduction}\label{sec1}

In various fields such as engineering \cite{R_m2013}, finance \cite{Z_m2015}, environmental analysis \cite{F_o2001}, management science \cite{T_a2010} and machine learning \cite{J_m2006}, one is faced with the problem that multiple objective have to be minimized concurrently, leading to a multiobjective optimization problem. These objective functions often conflict with each other, making it impossible to find a unique solution that optimizes all the objectives simultaneously. Instead, there exists a set of solutions known as the Pareto optimal solutions, which are characterized by the fact that an improvement in one objective will result in a deterioration in at least one of the other objectives.

One of the earliest and most widely used techniques for dealing with multiobjective optimization problems is the so-called scalarization method (see, e.g., \cite{miettinen1999,eichfelder2009}), which transforms the original problem into a single objective optimization problem. The main drawback of scalarization methods is the necessity of defining appropriate parameters to obtain a ``good'' scalarizing function. This requires a deep understanding of the problem's structure, which may not always be readily accessible. A class of parameter-free multiobjective optimization algorithms, known as descent algorithms, has garnered significant attention for not relying on prior information about the objective functions. This class of algorithms includes the steepest descent method \cite{fliege2000steepest,drummond2005steepest}, the Newton method \cite{fliege2009newton,wang2019extended}, the quasi-Newton method \cite{povalej2014quasi,lapucci2023limited,qu2011quasi,prudente2022quasi,ansary2015modified}, the trust-region method \cite{carrizo2016trust,ramirez2022}, the Barzilai-Borwein method \cite{morovati2016barzilai}, the conjugate gradient method \cite{lucambio2018nonlinear,goncalves2020extension,goncalves2022study}, the conditional gradient method \cite{AFP2021,chen2023conditional}, the proximal gradient method \cite{tanabe2019proximal,tanabe2023an}, and others.

In this paper, we are concerned with the following multiobjective optimization problem:
\begin{equation}\label{mop}
	\min_{x\in\mathbb{R}^{n}}\quad F(x)=(f_{1}(x),f_{2}(x),...,f_{m}(x))^{\top},\\
\end{equation}
where $F:\mathbb{R}^{n}\rightarrow\mathbb{R}^{m}$ is a continuously differentiable vector-valued function. As mentioned in \cite{upadhyay2024quasi}, a class of interval-valued multiobjective optimization problems can also be reformulated as the problem \eqref{mop}. In the present work, we focus on the conjugate gradient (CG) methods for solving \eqref{mop}. It is noteworthy that the CG methods for multiobjective optimization, initially proposed by Lucambio P\'{e}rez and Prudente \cite{lucambio2018nonlinear}, generates a sequence of iterates according to the following form:
\begin{equation}\label{iteration}
	x^{k+1}=x^{k}+t_{k}d^{k},\quad k=0,1,2,\ldots,
\end{equation}
where $d^{k}\in\mathbb{R}^{n}$ is the search direction, and $t_{k}>0$ is the stepsize. The direction in \cite{lucambio2018nonlinear} is defined by 
\begin{equation}\label{dk}
	d^{k}=\left\{\begin{array}{lllll}
		v(x^{k}), &\quad \text{if}~~k=0,\\
		v(x^{k})+\beta_{k}d^{k-1}, &\quad \text{if}~~k\geq 1,
	\end{array}\right.
\end{equation}
where $v(x^{k})$ is the multiobjective steepest descent direction (see \cite{fliege2000steepest} or \eqref{opt_sol}), and $\beta_{k}$ is a scalar algorithmic parameter determined by different formulas, each corresponding to a specific CG method. Lucambio P\'{e}rez and Prudente \cite{lucambio2018nonlinear} introduced the vector extensions of Fletcher--Reeves (FR), Conjugate descent (CD), Dai--Yuan (DY), Polak--Ribi\`{e}re--Polyak (PRP) and Hestenes--Stiefel (HS), which take the forms respectively defined by
\begin{equation*}
	\begin{aligned}
		\beta_{k}^{\rm FR}&=\dfrac{\psi(x^{k},v(x^{k}))}{\psi(x^{k-1},v(x^{k-1}))}, 	\\
		\beta_{k}^{\rm CD}&=\dfrac{\psi(x^{k},v(x^{k}))}{\psi(x^{k-1},d^{k-1})},\\
		\beta_{k}^{\rm DY}&=\dfrac{-\psi(x^{k},v(x^{k}))}{\psi(x^{k},d^{k-1})-\psi(x^{k-1},d^{k-1})},\\
		\beta_{k}^{\rm PRP}&=\dfrac{-\psi(x^{k},v(x^{k}))+\psi(x^{k-1},v(x^{k}))}{-\psi(x^{k-1},v(x^{k-1}))},\\
		\beta_{k}^{\rm HS}&=\dfrac{-\psi(x^{k},v(x^{k}))+\psi(x^{k-1},v(x^{k}))}{\psi(x^{k},d^{k-1})-\psi(x^{k-1},d^{k-1})},
	\end{aligned}
\end{equation*}
where $\psi(\cdot,\cdot)$ is defined in the next section. As reported in \cite{lucambio2018nonlinear}, in all cases, the vector extensions retrieve the classical parameters in the scalar minimization case. 
For the ways of selecting the parameter $\beta_{k}$, there are also some interesting works, such as the vector extension of Hager--Zhang (see \cite{goncalves2022study}), the vector extension of Liu-Storey (see \cite{goncalves2020extension}), the variants of FR, CD and PRP (see \cite{yahaya2023hybrid}), the alternative extension of HZ (see \cite{hu2024alternative}) and the vector extension of Dai-Liao (see \cite{yahaya2024dai}). Very recently, the authors in \cite{najafi2023multiobjective,chen2023nonlinear} have presented the CG methods for multiobjective optimization on Riemannian manifolds.

The global convergence of CG methods in multiobjective optimization crucially depends on the selection of the stepsize $t_{k}$. In \cite{lucambio2018nonlinear,goncalves2020extension,goncalves2022study,yahaya2023hybrid,hu2024alternative,yahaya2024dai}, it is usually required that the stepsize $t_{k}$ satisfies the standard Wolfe conditions
\begin{align}
	F(x^{k}+t_{k} d^{k})&\preceq F(x^{k})+\rho_{1} t_{k} \psi(x^{k},d^{k}),\label{armijo}\\
	\psi(x^{k}+t_{k}d^{k},d^{k})&\geq\rho_{2}\psi(x^{k},d^{k})\notag,
\end{align}
where ``$\preceq$'' is a partial order relation defined in the next section and $0<\rho_{1}<\rho_{2}<1$, or the strong Wolfe conditions, i.e., \eqref{armijo} and
\begin{equation}\label{wolfeb}
	\begin{aligned}
		\lvert\psi(x^{k}+t_{k}d^{k},d^{k})\rvert\leq\rho_{2}\lvert\psi(x^{k},d^{k})\rvert.
	\end{aligned}
\end{equation}
However, the implementation of such Wolfe-type stepsize strategies gives rise to certain potential concerns. Firstly, it is noteworthy that, in \eqref{armijo}, $m$ inequalities must be satisfied simultaneously, which becomes stricter as the number of objective functions increases (see \cite{morovati2016barzilai,mita2019nonmonotone}). Secondly, these types of line searches may involve  extensive computations in evaluating the value of objective functions, which can pose a substantial burden in certain problem scenarios. Additionally, the converging line search strategies are usually not easy to implement in practice. Recognizing these issues, it is advisable to steer clear of the line search procedures in the design of descent algorithms for multiobjective optimization.

In this paper, we introduce a novel stepsize strategy for multiobjective optimization, which encompasses the stepsize rule from our previous work \cite{chen2023memory}. The newly proposed stepsize scheme replaces the line search procedure with a fixed stepsize formula, motivated by the quadratic approximations of the objective functions, incorporating a sequence of positive definite matrices ${B^{k}}$. Here and below, the statement ``fixed
stepsize'' indicates ``without line search''. We establish the Zoutendijk-type condition for multiobjective optimization when using the proposed fixed
stepsize formula. This condition is a crucial component in deriving convergence results. We apply the proposed stepsize strategy to the CG methods \cite{lucambio2018nonlinear} with FR, CD, DY, PRP and HS parameters, and then establish the global convergence of these CG methods under certain assumptions. Finally, we present numerical results for these CG methods without line search and compare them with their counterparts employing stepsizes that satisfy the Wolfe-type conditions. 
It is noteworthy that, in the numerical experiments, we consider $B^{k}$ as a matrix with a simple structure that encapsulates some approximated  second-order information about the objective functions.

The remainder of this paper is as follows. In the next section, we present some basic definitions, notations and auxiliary results in multiobjective optimization. Section \ref{sec:3} gives the stepsize formula and its properties. Section \ref{sec:4} includes the main convergence results of the CG methods with the proposed stepsize strategy. Section \ref{sec:5} contains numerical experiments illustrating the performance of the proposed stepsize. Finally, in section \ref{sec:6}, some conclusions and remarks are given.

\section{Preliminaries}\label{sec:2}

Throughout this paper, we denote by $\langle \cdot,\cdot\rangle$ and $\|\cdot\|$, respectively, the usual inner product and the Euclidean norm in $\mathbb{R}^{n}$. For any positive integer $m$, let $\langle m\rangle=\{1,2,\ldots,m\}$. For any $d\in\mathbb{R}^{n}$, define $\|d\|_{A}=\sqrt{d^{\top}Ad}$, where $A$ is a positive definite matrix. Given a matrix $B\in\mathbb{R}^{m\times n}$, the norm of $B$ (see \cite{fliege2000steepest,hu2024alternative}) is computed as
\begin{equation*}
	\|B\|_{2,\infty}=\max_{x\neq0}\dfrac{\|Bx\|_{\infty}}{\|x\|}=\max_{i\in\langle m\rangle}\|B_{i,\cdot}\|=\max_{i\in\langle m\rangle}\left(\sum_{j=1}^{n}B_{i,j}^{2}\right)^{1/2}.
\end{equation*}
Let $\mathbb{R}_{+}^{m}$ and $\mathbb{R}_{++}^{m}$ denote the non-negative orthant and positive orthant of $\mathbb{R}^{m}$, respectively. As usual, for $u,v\in\mathbb{R}^{m}$, we use ``$\preceq$'' to denote the classical partial order defined by
\begin{equation*}
	u\preceq v ~\Leftrightarrow ~v-u\in\mathbb{R}_{+}^{n}.
\end{equation*}

\begin{definition}\normalfont\cite{miettinen1999}
	A point $\bar{x}\in\mathbb{R}^{n}$ is said to be \emph{Pareto optimal} of \eqref{mop} if there exists no $x\in\mathbb{R}^{n}$ such that $F(x)\preceq F(\bar{x})$ and $F(x)\neq F(\bar{x})$.
\end{definition}

A first-order necessary condition introduced in \cite{fliege2000steepest} for the Pareto optimality of a point $\bar{x}\in\mathbb{R}^{n}$ is
\begin{equation}\label{fir_ord_opt}
	JF(\bar{x})(\mathbb{R}^{n})\cap(-\mathbb{R}^{m}_{++})=\emptyset,
\end{equation}
where $JF(\bar{x})$ is the Jacobian (or the first-order derivative) of $F$ at $\bar{x}$, $ JF(\bar{x})(\mathbb{R}^{n})=\{JF(\bar{x})d:d\in\mathbb{R}^{n}\}$ and 
$JF(\bar{x})d=(\langle\nabla f_{1}(\bar{x}),d\rangle,\langle\nabla f_{2}(\bar{x}),d\rangle,\ldots,\langle\nabla f_{m}(\bar{x}),d\rangle)^{\top}.$

A point $\bar{x}\in\mathbb{R}^{n}$ satisfying the above relation \eqref{fir_ord_opt} is said to be \emph{Pareto critical} (see \cite{fliege2000steepest}).
Equivalently, for any $d\in\mathbb{R}^{n}$, there exists $i^{*}\in\langle m\rangle$ such that $(JF(\bar{x})d)_{i^{*}}=\langle\nabla f_{i^{*}}(\bar{x}),d\rangle\geq0,$
which implies $\max_{i\in\langle m\rangle}\langle\nabla f_{i}(\bar{x}),d\rangle\geq0$ for any $d\in\mathbb{R}^{n}$. 
Clearly, if $x\in\mathbb{R}^{n}$ is not a Pareto critical point, then there exists a vector $d\in\mathbb{R}^{n}$ satisfying $JF(x)d\in-\mathbb{R}^{m}_{++}$. We call the vector $d$ a \emph{descent direction} for $F$ at $x$.

Now define the function $\psi:\mathbb{R}^{n}\times\mathbb{R}^{n}\rightarrow \mathbb{R}$ by
\begin{equation}\label{h}
	\psi(x,d)=\max_{i\in\langle m\rangle}\langle\nabla f_{i}(x),d\rangle.
\end{equation}
Clearly, the function $\psi$ gives the following characterizations:
\begin{itemize}
	\item[-] $d\in\mathbb{R}^{n}$ is descent direction for $F$ at $x\in\mathbb{R}^{n}$ if and only if $\psi(x,d)<0$,
	\item[-] $x\in\mathbb{R}^{n}$ is Pareto critical if and only if $\psi(x,d)\geq0$ for all $d\in\mathbb{R}^{n}$.
\end{itemize}

Let us now consider the following scalar optimization problem:
\begin{equation}\label{sub_pro}
	\min_{d\in\mathbb{R}^{n}}\:\psi(x,d)+\frac{1}{2}\|d\|^{2}.
\end{equation}
Obviously, the objective function in \eqref{sub_pro} is proper, closed and strongly convex. Therefore, \eqref{sub_pro} admits a unique optimal solution, referred to as the \emph{steepest descent direction} (see \cite{fliege2000steepest}). Denote the optimal solution of (\ref{sub_pro}) by $v(x)$, i.e.,
\begin{equation}\label{opt_sol}
	v(x)=\mathop{\rm argmin}_{d\in\mathbb{R}^{n}}\: \psi(x,d)+\frac{1}{2}\|d\|^{2},
\end{equation}
and let the optimal value of (\ref{sub_pro}) be defined as $\theta(x)$, i.e., 
\begin{equation}\label{opt_val}
	\theta(x)=\psi(x,v(x))+\frac{1}{2}\|v(x)\|^{2}.
\end{equation}


Let us now give a characterization of Pareto critical points of  \eqref{mop}, which will be used in our subsequent analysis.

\begin{proposition}{\rm\cite{fliege2000steepest}}\label{pa_sta_equ}
	Let $v(\cdot)$ and $\theta(\cdot)$ be as in \eqref{opt_sol} and \eqref{opt_val}, respectively. The following statements hold:
	\begin{enumerate}
		\renewcommand{\theenumi}{\roman{enumi}}
		\renewcommand{\labelenumi}{\rm(\theenumi)}
		\item if $x$ is a Pareto critical point of \eqref{mop}, then $v(x)= 0$ and $\theta(x)=0$;\label{pa_sta_equ1}
		\item if $x$ is not a Pareto critical point of \eqref{mop}, then $v(x)\neq 0$, $\theta(x)<0$ and $\psi(x,v(x))<-\|v(x)\|^{2}/2<0$;\label{pa_sta_equ2}
		\item $v(\cdot)$ is continuous.\label{pa_sta_equ3}
	\end{enumerate}
\end{proposition}

We finish the section with the following properties.

\begin{lemma}\label{h_property}{\rm\cite{fukuda2014survey}}
	For all $x,y\in\mathbb{R}^{n}$, $\varrho> 0$ and $b_{1},b_{2}\in\mathbb{R}^{n}$, we obtain
	\begin{enumerate}
		\renewcommand{\theenumi}{\roman{enumi}}
		\renewcommand{\labelenumi}{\rm(\theenumi)}
		\item $\psi(x,\varrho b_{1})=\varrho \psi(x,b_{1})$;\label{h_property1}
		\item $\psi(x,b_{1}+b_{2})\leq \psi(x,b_{1})+\psi(x,b_{2})$ and $\psi(x,b_{1})-\psi(x,b_{2})\leq\psi(x,b_{1}-b_{2}) $;\label{h_property2}
	\end{enumerate}
\end{lemma}

\begin{lemma}\label{lemma_1}{\rm\cite{lucambio2018nonlinear}}
	For any scalars $a$, $b$ and $\alpha \neq 0$, we have
	\begin{enumerate}
		\renewcommand{\theenumi}{\roman{enumi}}
		\renewcommand{\labelenumi}{\rm(\theenumi)}
		\item	$(a+b)^{2}\leq2a^{2}+2b^{2}$;\label{lemma_11}
		\item $(a+b)^{2}\leq(1+2\alpha^{2})a^{2}+(1+1/(2\alpha^{2}))b^{2}$.\label{lemma_12}
	\end{enumerate}
\end{lemma}

\section{New stepsize and its properties}\label{sec:3}

This section first presents a novel stepsize formula and then lists some relations, all of which are independent of the choice of the parameter $\beta_{k}$. Before describing these works, we make use of the following assumptions:

\begin{assumption}\label{A1}\normalfont
	The level set $\mathcal{L}=\{x\in\mathbb{R}^{n}:F(x)\preceq F(x^{0})\}$ is bounded, where $x^{0}\in\mathbb{R}^{n}$ is a given starting point.
\end{assumption}

\begin{assumption}\label{A2}\normalfont
	The Jacobian $JF$ is $L$-Lipschitz continuous on an open set $\mathcal{B}$ containing $\mathcal{L}$, i.e.,	$\|JF(x)- JF(y)\|_{2,\infty}\leq L\|x-y\|$
	for all $x,y\in\mathcal{B}$.
\end{assumption}

\begin{assumption}\label{A3}\normalfont
	Each function $f_{i}$ is strongly convex with constant $\mu_{i}>0$, i.e., 
	$
		(\nabla f_{i}(x)- \nabla f_{i}(y))^{\top}(x-y)\geq\mu_{i}\| x-y\|^{2}
	$
	for all $x,y\in\mathcal{B}$ and $i\in\langle m\rangle$.
\end{assumption}

\begin{remark}\normalfont
	The considered Assumption {\normalfont\ref{A1}} and Assumption \ref{A2} are natural extensions of those made for the scalar case, which are used in \cite{fliege2000steepest,lapucci2023limited,lucambio2018nonlinear,goncalves2020extension,goncalves2022study,yahaya2023hybrid,hu2024alternative,chen2023memory}. 
	Moreover, Assumption \ref{A2} implies
	$
		\|\nabla f_{i}(x)- \nabla f_{i}(y)\|\leq L\|x-y\|
	$
	for all $x,y\in\mathcal{B}$ and $i\in\langle m\rangle$. Indeed, 
	\begin{equation*}
		\|JF(x)- JF(y)\|_{2,\infty}=\max_{i\in\langle m\rangle}\|(JF(x)-JF(y))_{i,\cdot}\|=\max_{i\in\langle m\rangle}\|\nabla f_{i}(x)-\nabla f_{i}(y)\|.
	\end{equation*}
	Assumption {\normalfont\ref{A3}} is also a commonly used tool for obtaining convergence analysis in the multiobjective optimization literature (see \cite{fliege2009newton,fukuda2014survey,zeng2019rate,tanabe2023rate,fliege2019complexity}). Clearly, Assumption {\normalfont\ref{A3}} implies Assumption {\normalfont\ref{A1}}.
\end{remark}

According to Assumption \ref{A3}, we have the following proposition.

\begin{proposition}
	Let $\mu=\min_{i\in\langle m\rangle}\mu_{i}$.
	Then, Assumption {\normalfont\ref{A3}} implies the following relation:
	\begin{equation}\label{str_con_pro}
		\psi(x^{k+1},x^{k+1}-x^{k})-\psi(x^{k},x^{k+1}-x^{k})\geq\mu\|x^{k+1}-x^{k}\|^{2}.
	\end{equation}
\end{proposition}
\noindent{\bf Proof.} Let $x=x^{k+1}$ and $y=x^{k}$ in Assumption \ref{A3}. For each $i\in\langle m\rangle$, we have
\begin{equation*}
	\begin{aligned}
		\nabla f_{i}(x^{k+1})^{\top}(x^{k+1}-x^{k})&\geq\nabla f_{i}(x^{k})^{\top}(x^{k+1}-x^{k})+\mu_{i}\| x^{k+1}-x^{k}\|^{2}\\
		&\geq\nabla f_{i}(x^{k})^{\top}(x^{k+1}-x^{k})+\mu\| x^{k+1}-x^{k}\|^{2}.
	\end{aligned}
\end{equation*}
Taking the $\max$ operator for all $i\in\langle m\rangle$ in the above inequality, and observing the definition of $\psi$, we obtain $$\psi(x^{k+1},x^{k+1}-x^{k})\geq\psi(x^{k},x^{k+1}-x^{k})+\mu\|x^{k+1}-x^{k}\|^{2},$$
which concludes the proof.
\qed

Now, we introduce the fixed stepsize strategy for multiobjective optimization. For each $i\in\langle m\rangle$, we consider the following quadratic model: 
\begin{equation}\label{quadratic}
	\begin{aligned}
		f_{i}(x^{k}+t d^{k})-f_{i}(x^{k})\approx t\nabla f_{i}(x^{k})^{\top}d^{k}+\dfrac{t^{2}}{2}\|d^{k}\|^{2}_{B^{k}},
	\end{aligned}
\end{equation}
where $B^{k}$ is an $n\times n$ positive definite symmetric matrix. From the definition of $\psi$ in \eqref{h}, the right side of \eqref{quadratic} is less than or equal to $t\psi(x^{k},d^{k})+t^{2}\|d^{k}\|^{2}_{B^{k}}/2$, which is a convex quadratic function that has the optimal solution
$t=-\psi(x^{k},d^{k})/\|d^{k}\|^{2}_{B^{k}}.$
Upon recognizing this issue, we determine the stepsize in this work by using the following formula rather than the line search procedure:
\begin{equation}\label{tk}
	t_{k}=-\frac{\delta\psi(x^{k},d^{k})}{\|d^{k}\|^{2}_{B^{k}}},
\end{equation}
where $\delta>0$ is a scalar. Clearly, if $d^{k}$ is a descent direction, then $\psi(x^{k},d^{k})<0$, and thus $t_{k}>0$.

\begin{remark}\normalfont
	If $B^{k}=I$  (the identity matrix) for all $k$ and $\delta=1/L$, then \eqref{tk} reduces to the stepsize-II strategy provided in \cite{chen2023memory}. If $m=1$,  we fall back to single objective optimization, and then \eqref{tk} retrieves to the stepsize strategy proposed by Sun and Zhang \cite{sun2001global}, which has been further investigated in \cite{chen2002global,zhang2005a,zhu2008global,sendilkkumaar2018}.
\end{remark}

If the stepsize is given by formula \eqref{tk}, then the following conclusions hold.

\begin{lemma}\label{lemma1}
	Consider an iteration of the form \eqref{iteration}, where $d^{k}$ is a descent direction for $F$ at $x^{k}$ and $t_{k}$ satisfies \eqref{tk}. For all $k\geq 0$, we have	
	\begin{equation}\label{lemma1_1}
		\psi(x^{k+1},d^{k})=\rho_{k}\psi(x^{k},d^{k}),
	\end{equation}
	where 
	\begin{equation}\label{rho_k}
		\rho_{k}=1-\dfrac{\delta\eta_{k}\|d^{k}\|^{2}}{\|d^{k}\|_{B^{k}}^{2}}
	\end{equation}
	and 
	\begin{equation*}
		\eta_{k}=\left\{\begin{array}{lllll}
			0, &\quad \text{if}~~t_{k}=0,\\
			\dfrac{\psi(x^{k+1},x^{k+1}-x^{k})-\psi(x^{k},x^{k+1}-x^{k})}{\|x^{k+1}-x^{k}\|^{2}}, &\quad \text{if}~~ t_{k}\neq 0.
		\end{array}\right.
	\end{equation*}
\end{lemma}
\noindent{\bf Proof.} It is evident that \eqref{lemma1_1} holds when $t_{k}=0$. Now consider the case of $t_{k}\neq 0$. By \eqref{iteration}, Proposition \ref{h_property}\eqref{h_property1} and the definition of $t_{k}$ in \eqref{tk}, we have
\begin{equation}\label{lemma1_2}
	\begin{aligned}
		\psi(x^{k+1},d^{k})&=\psi(x^{k},d^{k})+ \psi(x^{k+1},d^{k})-\psi(x^{k},d^{k})\\
		&=\psi(x^{k},d^{k})+\frac{1}{t_{k}}(\psi(x^{k+1}, x^{k+1}-x^{k})-\psi(x^{k}, x^{k+1}-x^{k}))\\
		&=\psi(x^{k},d^{k})+\frac{\eta_{k}}{t_{k}}\|x^{k+1}-x^{k}\|^{2}\\
		&=\psi(x^{k},d^{k})+\eta_{k}t_{k}\|d^{k}\|^{2}\\
		&=\left(1-\dfrac{\delta\eta_{k}\|d^{k}\|^{2}}{\|d^{k}\|_{B^{k}}^{2}}\right)\psi(x^{k},d^{k}).
	\end{aligned}
\end{equation}
Therefore, if we set $\rho_{k}=1-\delta\eta_{k}\|d^{k}\|^{2}/\|d^{k}\|_{B^{k}}^{2}$ in \eqref{lemma1_2}, then \eqref{lemma1_1} is satisfied and the proof is complete.
\qed

\begin{assumption}\normalfont\label{A4}
	There exist positive constants $a_{\min}$ and $a_{\max}$ such that
	\begin{equation}\label{amin_amax}
		a_{\min}\|d\|^2\leq d^{\top}B^{k}d\leq a_{\max}\|d\|^2
	\end{equation}
	for all $d\in\mathbb{R}^{n}$ and all $k\geq0$. 
\end{assumption}

In the subsequent analysis, the parameter $\delta>0$ is selected such that 
\begin{equation}\label{delta}
	\delta\in\left(0,\dfrac{a_{\min}}{L}\right).
\end{equation}
In this case, it can be seen that 
\begin{equation}\label{ldelta}
	\dfrac{L\delta}{a_{\min}}<1.
\end{equation}

\begin{corollary}\label{coro1}
	Let $\rho_{k}$ be as in \eqref{rho_k}. Then, the following statements hold:
	\begin{enumerate}
		\renewcommand{\theenumi}{\roman{enumi}}
\renewcommand{\labelenumi}{\rm(\theenumi)}
		\item 	if Assumption {\normalfont\ref{A2}} holds, then 
		\begin{equation*}
			\rho_{k}\in\left[1-\dfrac{L\delta}{a_{\min}},1+\dfrac{L\delta}{a_{\min}}\right];
		\end{equation*}\label{coro11}
		\item 	if Assumption {\normalfont\ref{A3}} holds, then 
		\begin{equation*}
			\rho_{k}\in\left(0,1-\dfrac{\mu\delta}{a_{\max}}\right].
		\end{equation*}\label{coro12}
	\end{enumerate}
\end{corollary}
\noindent{\bf Proof.}  
From the definition of $\rho_{k}$ in \eqref{rho_k}, it follows that
\begin{equation}\label{coro2}
	\begin{aligned}
		\rho_{k}
		=1-
		\dfrac{\delta(\psi(x^{k+1},x^{k+1}-x^{k})-\psi(x^{k},x^{k+1}-x^{k}))}{\|x^{k+1}-x^{k}\|^{2}}\dfrac{\|d^{k}\|^{2}}{\|d^{k}\|_{B^{k}}^{2}}.
	\end{aligned}
\end{equation}
By \cite[Lemma 4.1]{hu2024alternative} and Assumption \ref{A2}, we  have
\begin{equation*}
	\begin{aligned}
		\lvert\psi(x^{k+1},x^{k+1}-x^{k})-\psi(x^{k},x^{k+1}-x^{k})\rvert
		&\leq\|JF(x^{k+1})-JF(x^{k})\|_{2,\infty}\|x^{k+1}-x^{k}\|\\
		&\leq L\|x^{k+1}-x^{k}\|^{2},
	\end{aligned}
\end{equation*}
Therefore,
\begin{equation*}
	\dfrac{\psi(x^{k+1},x^{k+1}-x^{k})-\psi(x^{k},x^{k+1}-x^{k})}{\|x^{k+1}-x^{k}\|^{2}}\in[-L,L].
\end{equation*}
This together with \eqref{coro2} and \eqref{amin_amax} gives us the bounds for $\rho_{k}$ in item (i).

The item (ii) follows from \eqref{str_con_pro}, \eqref{amin_amax} and \eqref{coro2}.
\qed

Next, we establish that the iterative form satisfies a condition of Zoutendijk's type, which  is subsequently used to prove the global convergence of the CG methods with our stepsize strategy.

\begin{lemma}\label{zou1}
	Suppose that Assumptions {\normalfont\ref{A1}}, {\normalfont\ref{A2}} and {\normalfont\ref{A4}} hold. Consider an iteration of the form \eqref{iteration}, where $d^{k}$ is a descent direction for $F$ at $x^{k}$ and $t_{k}$ satisfies \eqref{tk}. Then
	\begin{equation}\label{zou11}
		\sum_{k\geq0}\frac{\psi^{2}(x^{k},d^{k})}{\|d^{k}\|^{2}}<\infty.
	\end{equation}
\end{lemma}
\noindent{\bf Proof.}  
It follows from the mean-value theorem that 
\begin{equation}\label{zou2}
	f_{i}(x^{k+1})-f_{i}(x^{k})=\nabla f_{i}(\xi_{i}^{k})^{\top}(x^{k+1}-x^{k})
\end{equation}
for each $i\in\langle m\rangle$, where $\xi_{i}^{k}\in\mathbb{R}^{n}$ lies in the line segment connecting $x^{k}$ and $x^{k+1}$. By the Cauchy-Schwarz inequality, Assumption \ref{A2}, the definition of $\psi$, Lemma \ref{h_property}\eqref{h_property1}, \eqref{tk} and \eqref{amin_amax}, for each $i\in\langle m\rangle$, we have
\begin{equation*}
	\begin{aligned}
		\nabla f_{i}(\xi_{i}^{k})^{\top}(x^{k+1}-x^{k})&=\nabla f_{i}(x^{k})^{\top}(x^{k+1}-x^{k})+(\nabla f_{i}(\xi_{i}^{k})-\nabla f_{i}(x^{k}))^{\top}(x^{k+1}-x^{k})\\
		&\leq\nabla f_{i}(x^{k})^{\top}(x^{k+1}-x^{k})+\|\nabla f_{i}(\xi_{i}^{k})-\nabla f_{i}(x^{k})\|\|x^{k+1}-x^{k}\|\\
		&\leq\nabla f_{i}(x^{k})^{\top}(x^{k+1}-x^{k})+L\|x^{k+1}-x^{k}\|^{2}\\
		&\leq t_{k}\psi(x^{k},d^{k})+Lt_{k}^{2}\|d^{k}\|^{2}\\
		&=-t_{k}\psi(x^{k},d^{k})\left(\dfrac{L\delta\|d^{k}\|^{2}}{\|d^{k}\|_{B^{k}}^{2}}-1\right)\\
		&= \dfrac{\delta\psi^{2}(x^{k},d^{k})}{\|d^{k}\|_{B^{k}}^{2}}\left(\dfrac{L\delta\|d^{k}\|^{2}}{\|d^{k}\|_{B^{k}}^{2}}-1\right)\\
		&\leq\left(\dfrac{L\delta^{2}}{a_{\min}}-\delta\right)\dfrac{\psi^{2}(x^{k},d^{k})}{\|d^{k}\|_{B^{k}}^{2}},
	\end{aligned}
\end{equation*}
which, combined with \eqref{zou2}, yields
\begin{equation}\label{F_descent}
	f_{i}(x^{k+1})-f_{i}(x^{k})\leq\delta\left(\dfrac{L\delta}{a_{\min}}-1\right)\dfrac{\psi^{2}(x^{k},d^{k})}{\|d^{k}\|_{B^{k}}^{2}},
\end{equation}
for each $i\in\langle m\rangle$. This implies that $\{F(x^{k})\}_{k\geq 0}$ is monotone decreasing because $L\delta/a_{\min}< 1$ in \eqref{ldelta} and that
\begin{equation*}
	\sum_{j=0}^{k}\dfrac{\psi^{2}(x^{j},d^{j})}{\|d^{j}\|_{B^{j}}^{2}}\leq \dfrac{a_{\min}}{a_{\min}\delta-L\delta^{2}}(f_{i}(x^{0})-f_{i}(x^{k+1}))
\end{equation*}
for each $i\in\langle m\rangle$ and all $k\geq0$. Then, by Assumption {\normalfont\ref{A1}} and the continuity arguments, since $\{x^{k}\}\subset\mathcal{L}$, there exists $\bar{F}\in\mathbb{R}^{m}$ such that $\bar{F}\preceq F(x^{k})$ for all $k\geq0$ and 
\begin{equation}\label{zou3}
	\sum_{j=0}^{k}\dfrac{\psi^{2}(x^{j},d^{j})}{\|d^{j}\|_{B^{j}}^{2}}\leq \dfrac{a_{\min}}{a_{\min}\delta-L\delta^{2}}(f_{i}(x^{0})-\bar{F}_{i})
\end{equation}
for each $i\in\langle m\rangle$  and all $k\geq0$. Thus, by \eqref{amin_amax} and \eqref{zou3}, we obtain
\begin{equation*}
	\sum_{k\geq0}\dfrac{\psi^{2}(x^{k},d^{k})}{\|d^{k}\|^{2}}\leq a_{\max}\sum_{k\geq0}\dfrac{\psi^{2}(x^{k},d^{k})}{\|d^{k}\|_{B^{k}}^{2}}<\infty,
\end{equation*}
which concludes the proof.
\qed

\begin{remark}\normalfont
	In \cite[Proposition 3.3]{lucambio2018nonlinear}, the condition of Zoutendijk's type was establised using the standard Wolfe conditions. In contrast, our result does not require these conditions. If we replace Assumption {\normalfont\ref{A1}} with Assumption {\normalfont\ref{A3}} in the above lemma, the corresponding conclusion still hold, because Assumption {\normalfont\ref{A3}} implies Assumption {\normalfont\ref{A1}}.
\end{remark}

\begin{remark}\label{re4.1}\normalfont
	The multiobjective steepest descent method\footnote{In this method, the Armijo line search \cite{fliege2000steepest} is proposed to ensure the convergence of the sequence of iterates generated by the method. Specifically, the stepsize $t_{k}$ is chosen as the largest value from $\{(1/2)^{0},(1/2)^{1},\ldots\}$ such that
	$F(x^{k}+t_{k} d^{k})\preceq F(x^{k})+\beta t_{k} JF(x^{k})d^{k}$.
	} proposed by Fliege and Svaiter	\cite{fliege2000steepest} converges in the sense that $\lim_{k\rightarrow\infty}\|v(x^{k})\|=0$, provided that it employs the stepsize rule \eqref{tk} instead of the Armijo line search. Indeed, consider the iteration \eqref{iteration} with $d^{k}=v(x^{k})$ and assume that $t_{k}$ satisfies \eqref{tk} for all $k\geq0$. According to the relation $-\psi(x^{k},v(x^{k}))\geq\|v(x^{k})\|^{2}/2$ in Proposition \ref{pa_sta_equ}\eqref{pa_sta_equ2}, we have $$\dfrac{\psi^{2}(x^{k},d^{k})}{\|d^{k}\|^{2}}=\dfrac{\psi^{2}(x^{k},v(x^{k}))}{\|v(x^{k})\|^{2}}\geq\dfrac{\|v(x^{k})\|^{2}}{4}.$$
	Summing the above inequality over all $k$, and combining it with \eqref{zou11}, we have $\sum_{k\geq0}\|v(x^{k})\|^{2}<\infty$, and thus $\|v(x^{k})\|\rightarrow0$ as $k\rightarrow\infty$.
\end{remark}
%

\section{Convergence analysis}\label{sec:4}

In this section, we analyze the convergence properties of the CG methods with the vector extensions of FR, CD, DY, PRP and HS parameters and the proposed stepsize rule. We now present the CG algorithm for the multiobjective optimization problem \eqref{mop}.
\begin{algorithm}
	\caption{CG Algorithm to solve \eqref{mop}}\label{alg1}
	\begin{algorithmic}[1]  
		\State \textbf{Input:} $x^0 \in \mathbb{R}^n$.
		\For{$k=0,1,2,\ldots$}
			\State Compute $v(x^k)$;
			\If{$v(x^k) = 0$}
				\State \textbf{Return}
				Pareto critical point $x^{k}$;
			\EndIf
			\State Select a parameter $\beta_{k}$;
			\State Compute $d^k$ using \eqref{dk};
			\State Compute $t_k$ using \eqref{tk};
			\State Update $x^{k+1} = x^k + t_k d^k$;
		\EndFor
	\end{algorithmic}
\end{algorithm}

\begin{lemma}\label{zout}
	Consider Algorithm {\normalfont\ref{alg1}} and let Assumptions {\normalfont\ref{A1}}, {\normalfont\ref{A2}} and {\normalfont\ref{A4}} hold. If
	$\beta_{k}\geq0$, $d^{k}$ is a descent direction of $F$ at $x^{k}$, and $\lim\inf_{k\rightarrow\infty}\|v(x^{k})\|\neq 0$, then
	\begin{equation}\label{zout00}
		\sum_{d^{k}\neq0}\frac{\psi^{2}(x^{k},v(x^{k}))}{\|d^{k}\|^{2}}<\infty.
	\end{equation}
\end{lemma}
\noindent{\bf Proof.}  
It follows from $\lim\inf_{k\rightarrow\infty}\|v(x^{k})\|\neq 0$ that there exists $\gamma>0$ such that 
\begin{equation}\label{lem5}
	\|v(x^{k})\|\geq\gamma
\end{equation}
for all $k\geq0$. Set $\tau_{k}=\lvert\psi(x^{k},d^{k})\rvert/\|d^{k}\|$. From Lemma \ref{zou1}, we have $\tau_{k}\leq\gamma/8$ for all large $k$ and
\begin{equation}\label{zout0}
	\sum_{k\geq0}\tau_{k}^{2}<\infty.
\end{equation} 
From \eqref{lemma1_1}, Corollary \ref{coro1}\eqref{coro11} and \eqref{ldelta}, we obtain
\begin{equation}\label{zout1}
	\begin{aligned}
		\lvert\psi(x^{k},d^{k-1})\rvert&=\lvert\rho_{k-1}\psi(x^{k-1},d^{k-1})\rvert\\
		&\leq\left(1+\dfrac{L\delta}{a_{\min}}\right)\lvert\psi(x^{k-1},d^{k-1})\rvert\\
		&<2\lvert\psi(x^{k-1},d^{k-1})\rvert.
	\end{aligned}
\end{equation}
By the definition of $d^{k}$ and Lemma \ref{h_property}\eqref{h_property1}--\eqref{h_property2}, one has
\begin{equation*}
	\begin{aligned}
		\psi(x^{k},d^{k})=\psi(x^{k},v(x^{k})+\beta_{k}d^{k-1})\leq\psi(x^{k},v(x^{k}))+\beta_{k}\psi(x^{k},d^{k-1}),
	\end{aligned}
\end{equation*}
and thus
\begin{equation}\label{zout2}
	\begin{aligned}
		\dfrac{-\psi(x^{k},v(x^{k}))}{\|d^{k}\|}
		&\leq\dfrac{-\psi(x^{k},d^{k})+\beta_{k}\psi(x^{k},d^{k-1})}{\|d^{k}\|}\\
		&\leq\dfrac{\lvert\psi(x^{k},d^{k})\rvert}{\|d^{k}\|}+\dfrac{\beta_{k}\lvert\psi(x^{k},d^{k-1})\rvert}{\|d^{k}\|}\\
		&\leq\tau_{k}+\dfrac{2\beta_{k}\lvert\psi(x^{k-1},d^{k-1})\rvert}{\|d^{k}\|}\\
		&=\tau_{k}+2\tau_{k-1}\dfrac{\|\beta_{k}d^{k-1}\|}{\|d^{k}\|}\\
		&\leq\tau_{k}+2\tau_{k-1}\dfrac{\|d^{k}\|+\|v(x^{k})\|}{\|d^{k}\|}\\
		&=\tau_{k}+2\tau_{k-1}+2\tau_{k-1}\dfrac{\|v(x^{k})\|^{2}}{\|d^{k}\|\|v(x^{k})\|}\\
		&\leq \tau_{k}+2\tau_{k-1}+2\tau_{k-1}\dfrac{\|v(x^{k})\|^{2}}{\gamma\|d^{k}\|}\\
		&\leq \tau_{k}+2\tau_{k-1}+2\left(\dfrac{\gamma}{8}\right)\dfrac{-2\psi(x^{k},v(x^{k}))}{\gamma\|d^{k}\|}\\
		&=\tau_{k}+2\tau_{k-1}+\dfrac{-\psi(x^{k},v(x^{k}))}{2\|d^{k}\|},
	\end{aligned}
\end{equation}
where the third inequality follows from \eqref{zout1}, the first equality is due to the definition of $\tau_{k}$, the fifth inequality follows from \eqref{lem5} and the sixth inequality holds in view of the relation $-\psi(x^{k},v(x^{k}))\geq\|v(x^{k})\|^{2}/2$ and the fact that $\tau_k\leq \gamma/8$. With a rearrangement of \eqref{zout2}, we can obtain
\begin{equation*}
	0\leq\dfrac{-\psi(x^{k},v(x^{k}))}{\|d^{k}\|}\leq2\tau_{k}+4\tau_{k-1}\leq4(\tau_{k}+\tau_{k-1}).
\end{equation*}
This, combined with Lemma \ref{lemma_1}\eqref{lemma_12} and \eqref{zout0}, yields
\begin{equation*}
	\sum_{d^{k}\neq0}\dfrac{\psi^{2}(x^{k},v(x^{k}))}{\|d^{k}\|^{2}}\leq\sum_{d^{k}\neq0}16(\tau_{k}+\tau_{k-1})^{2}\leq\sum_{d^{k}\neq0}32(\tau_{k}^{2}+\tau_{k-1}^{2})<\infty,
\end{equation*}
which concludes the proof.
\qed

\begin{remark}\normalfont
	If we substitute Assumption {\normalfont\ref{A1}} with Assumption {\normalfont\ref{A3}} in the lemma above, the conclusion remains valid since Assumption {\normalfont\ref{A3}} implies Assumption {\normalfont\ref{A1}}.
\end{remark}

The following three theorems respectively discuss the convergence results of Algorithm {\normalfont\ref{alg1}} with the fixed stepsize formula \eqref{tk} under the FR, CD and DY parameters.

\begin{theorem}\label{fr_theorem}
	Consider Algorithm {\normalfont\ref{alg1}} and let $0<\xi< 1$. Assume that Assumptions {\normalfont\ref{A1}}, {\normalfont\ref{A2}} and {\normalfont\ref{A4}} hold. If $\lvert\beta_{k}\rvert\leq\xi\beta_{k}^{\rm FR}$ and $d^{k}$ is a descent direction of $F$ at $x^{k}$, then $\lim\inf_{k\rightarrow\infty}\|v(x^{k})\|=0$.
\end{theorem}
\noindent{\bf Proof.}  
Assume on the contrary that there exists $\gamma>0$ such that
\begin{equation*}
	\|v(x^{k})\|\geq\gamma
\end{equation*}
for all $k\geq0$. 
By the definition of $d^{k}$, and Lemma \ref{lemma_1}\eqref{lemma_12} with $a=\|v(x^{k})\|$, $b=\lvert\beta_{k}\rvert\|d^{k-1}\|$ and $\alpha=\xi/\sqrt{2(1-\xi^{2})}$, one has
\begin{equation*}
	\begin{aligned}
		\|d^{k}\|^{2}&=\|v(x^{k})+\beta_{k}d^{k-1}\|^{2}\\
		&\leq(\|v(x^{k})\|+\lvert\beta_{k}\rvert\|d^{k-1}\|)^{2}\\
		&\leq\dfrac{1}{1-\xi^{2}}\|v(x^{k})\|^{2}+\dfrac{1}{\xi^{2}}\beta_{k}^{2}\|d^{k-1}\|^{2}.
	\end{aligned}
\end{equation*}
Dividing both sides by $\psi^{2}(x^{k},v(x^{k}))$ in the above relation, and observing that $$\lvert\beta_{k}\rvert\leq\xi\beta_{k}^{\rm FR}=\dfrac{\xi \psi(x^{k},v(x^{k}))}{\psi(x^{k-1},v(x^{k-1}))},$$ we have
\begin{equation*}
	\dfrac{\|d^{k}\|^{2}}{\psi^{2}(x^{k},v(x^{k}))}\leq\dfrac{1}{1-\xi^{2}}\frac{\|v(x^{k})\|^{2}}{\psi^{2}(x^{k},v(x^{k}))}+\dfrac{\|d^{k-1}\|^{2}}{\psi^{2}(x^{k-1},v(x^{k-1}))}.
\end{equation*}
This together with $\gamma^{2}\leq\|v(x^{k})\|^{2}\leq-2\psi(x^{k},v(x^{k}))$ gives us
\begin{equation*}
	\begin{aligned}
		\dfrac{\|d^{k}\|^{2}}{\psi^{2}(x^{k},v(x^{k}))}&\leq\dfrac{\|d^{k-1}\|^{2}}{\psi^{2}(x^{k-1},v(x^{k-1}))}+\frac{4}{\gamma^{2}(1-\xi^{2})}\\
		&\leq\dfrac{\|d^{k-2}\|^{2}}{\psi^{2}(x^{k-2},v(x^{k-2}))}+\dfrac{8}{\gamma^{2}(1-\xi^{2})}\\
		&\leq\cdots\\
		&\leq\dfrac{\|d^{0}\|^{2}}{\psi^{2}(x^{0},v(x^{0}))}+\dfrac{4k}{\gamma^{2}(1-\xi^{2})}\\
		&\leq\dfrac{ 4}{\gamma^{2}}+\dfrac{4k}{\gamma^{2}(1-\xi^{2})}.
	\end{aligned}
\end{equation*}
Thus, we obtain
\begin{equation*}
	\begin{aligned}
		\sum_{d^{k}\neq0}\dfrac{\psi^{2}(x^{k},v(x^{k}))}{\|d^{k}\|^{2}}\geq\sum_{d^{k}\neq0}\dfrac{\gamma^{2}(1-\xi^{2})}{4}\dfrac{1}{k+1}=\infty,
	\end{aligned}
\end{equation*}
contradicting the relation \eqref{zout00} and concluding the proof.
\qed

\begin{theorem}\label{cd_theorem}
	Consider Algorithm {\normalfont\ref{alg1}} with $\beta_{k}=\beta_{k}^{\rm CD}\geq0$. If Assumptions {\normalfont\ref{A1}}, {\normalfont\ref{A2}} and {\normalfont\ref{A4}} hold, then $\lim\inf_{k\rightarrow\infty}\|v(x^{k})\|=0$.
\end{theorem}
\noindent{\bf Proof.}  
According to the definitions of $d^{k}$ and $\beta_{k}^{\rm CD}$, Lemma \ref{h_property}\eqref{h_property1}--\eqref{h_property2}, \eqref{lemma1_1} and Corollary \ref{coro1}, one has
\begin{equation}\label{cd_them1}
	\begin{aligned}
		\psi(x^{k},d^{k})&=\psi\left(x^{k},v(x^{k})+\frac{\psi(x^{k},v(x^{k}))}{\psi(x^{k-1},d^{k-1})}d^{k-1}\right)\\
		&\leq\psi(x^{k},v(x^{k}))+\frac{\psi(x^{k},v(x^{k}))}{\psi(x^{k-1},d^{k-1})}\psi(x^{k},d^{k-1})\\
		&=\psi(x^{k},v(x^{k}))+\frac{\psi(x^{k},v(x^{k}))}{\psi(x^{k-1},d^{k-1})}(\rho_{k-1}\psi(x^{k-1},d^{k-1}))\\
		&=(1+\rho_{k-1})\psi(x^{k},v(x^{k}))
		<0,
	\end{aligned}
\end{equation}
where $\rho_{k-1}>0$, this means that $d^{k}$ is a descent direction. From \eqref{cd_them1} and Corollary \ref{coro1}, and observing $\psi(x^{k},v(x^{k}))<0$, it follows that
\begin{equation*}
	\begin{aligned}
		\beta_{k}^{\rm CD}&=\frac{\psi(x^{k},v(x^{k}))}{\psi(x^{k-1},d^{k-1})}\\
		&\leq\dfrac{\psi(x^{k},v(x^{k}))}{(1+\rho_{k-2})\psi(x^{k-1},v(x^{k-1}))}\\
		&=\dfrac{1}{1+\rho_{k-2}}\beta_{k}^{\rm FR}\\
		&\leq\dfrac{a_{\min}}{2a_{\min}-L\delta}\beta_{k}^{\rm FR}.
	\end{aligned}
\end{equation*}
Since $\delta\in(0,a_{\min}/L)$, we have $a_{\min}/(2a_{\min}-L\delta)\in(0,1).$ The rest of the proof follows from Theorem \ref{fr_theorem}.
\qed

\begin{theorem}\label{dy_theorem}
	Consider Algorithm {\normalfont\ref{alg1}} with $\beta_{k}=\beta_{k}^{\rm DY}\geq0$. Let Assumptions {\normalfont\ref{A2}}--{\normalfont\ref{A4}} hold and $c=1-\mu\delta/a_{\max}$. Then, $d^{k}$ is a descent direction. Furthermore, if $$\beta_{k}=\eta\beta_{k}^{\rm DY}, \quad where \quad 0<\eta<\dfrac{1-c}{1+c},$$ then $\lim\inf_{k\rightarrow\infty}\|v(x^{k})\|=0$.
\end{theorem}
\noindent{\bf Proof.}  
It follows from the definition of $\beta_{k}^{\rm DY}$ 
and \eqref{lemma1_1} that
\begin{equation}\label{dy_lemma_1}
	\beta_{k}^{\rm DY}=\dfrac{-\psi(x^{k},v(x^{k}))}{\psi(x^{k},d^{k-1})-\psi(x^{k-1},d^{k-1})}=\dfrac{\psi(x^{k},v(x^{k}))}{(1-\rho_{k-1})\psi(x^{k-1},d^{k-1})}.
\end{equation}
This, combined with the definition of $d^{k}$, Lemma \ref{h_property}\eqref{h_property1}--\eqref{h_property2} and \eqref{lemma1_1}, gives
\begin{equation*}\label{dy_lemma_1_1}
	\begin{aligned}
		\psi(x^{k},d^{k})&=\psi(x^{k},v(x^{k})+\beta_{k}^{\rm DY}d^{k-1})\\
		&\leq\psi(x^{k},v(x^{k}))+\beta_{k}^{\rm DY}\psi(x^{k},d^{k-1})\\
		&=\psi(x^{k},v(x^{k}))+\dfrac{\psi(x^{k},v(x^{k}))}{(1-\rho_{k-1})\psi(x^{k-1},d^{k-1})}(\rho_{k-1}\psi(x^{k-1},d^{k-1}))\\
		&=\dfrac{1}{1-\rho_{k-1}}\psi(x^{k},v(x^{k})).
	\end{aligned}
\end{equation*}
Using Corollary \ref{coro1}\eqref{coro12}, we have $\rho_{k}\in(-c,c]$ with $0<c\leq 1$. Therefore, one has 
\begin{equation}\label{des}
	\psi(x^{k},d^{k})\leq\dfrac{1}{1+c}\psi(x^{k},v(x^{k}))<0.
\end{equation}
This implies that $d^{k}$ is a descent direction. 

On the other hand, from \eqref{lemma1_1} and \eqref{des}, it follows that
\begin{equation}\label{dy_lemma_2}
	\begin{aligned}
		\psi(x^{k},d^{k-1})-\psi(x^{k-1},d^{k-1})
		&=(\rho_{k-1}-1)\psi(x^{k-1},d^{k-1})\\
		&\geq(c-1)\psi(x^{k-1},d^{k-1})\\
		&\geq\dfrac{c-1}{1+c}\psi(x^{k-1},v(x^{k-1}))
		>0.
	\end{aligned}
\end{equation}
Let $\xi=\eta(1+c)/(1-c)$. Then $\xi\in(0,1)$. By the definition of $\beta_{k}$ and \eqref{dy_lemma_2}, we have
\begin{equation*}
	\begin{aligned}
		\beta_{k}&=\eta\beta_{k}^{\rm DY}\\
		&=\dfrac{\xi(1-c)}{1+c}\dfrac{-\psi(x^{k},v(x^{k}))}{\psi(x^{k},d^{k-1})-\psi(x^{k-1},d^{k-1})}\\
		&=\dfrac{\xi(1-c)}{1+c}\left(\dfrac{-\psi(x^{k},v(x^{k}))}{-\psi(x^{k-1},v(x^{k-1}))}\right)\left(\dfrac{-\psi(x^{k-1},v(x^{k-1}))}{\psi(x^{k},d^{k-1})-\psi(x^{k-1},d^{k-1})}\right)\\
		&\leq\xi\dfrac{\psi(x^{k},v(x^{k}))}{\psi(x^{k-1},v(x^{k-1}))}\\
		&=\xi\beta_{k}^{\rm FR}.
	\end{aligned}
\end{equation*}
The remainder of the proof follows the same steps as those in the proof of Theorem \ref{fr_theorem}.
\qed

Now, we analyze the convergence of the CG algorithm with the proposed stepsize rule \eqref{tk} under the PRP and HS parameters.

\begin{theorem}\label{prp_theorem}
	Suppose that Assumptions {\normalfont\ref{A1}}, {\normalfont\ref{A2}} and {\normalfont\ref{A4}} hold. Consider Algorithm {\normalfont\ref{alg1}} with $\beta_{k}=\max\{\beta_{k}^{\rm PRP},0\}$. If $d^{k}$ is a descent direction, then $\lim\inf_{k\rightarrow\infty}\|v(x^{k})\|=0$.
\end{theorem}
\noindent{\bf Proof.}  
Assume by contradiction that there exists a positive constant $\gamma$ such that 
$$\|v(x^{k})\|\geq\gamma$$ for all $k\geq 0$. 
From the definition of $\beta_{k}^{\rm PRP}$ and \eqref{lemma1_1}, and observing that $ \psi(x^{k},v(x^{k}))\leq-\|v(x^{k})\|^{2}/2$,
we have
\begin{equation*}
	\begin{aligned}
		\lvert\beta_{k}^{\rm PRP}\rvert&=\left\lvert\dfrac{-\psi(x^{k},v(x^{k}))+\psi(x^{k-1},v(x^{k}))}{-\psi(x^{k-1},v(x^{k-1}))}\right\rvert\\
		&\leq2\dfrac{\lvert-\psi(x^{k},v(x^{k}))+\psi(x^{k-1},v(x^{k}))\rvert}{\|v(x^{k-1})\|^{2}}\\
		&\leq\dfrac{2}{\gamma^{2}}\lvert-\psi(x^{k},v(x^{k}))+\psi(x^{k-1},v(x^{k}))\rvert.
	\end{aligned}
\end{equation*}

Using the iteration \eqref{iteration}, the definition of $t_{k}$ in \eqref{tk} and \eqref{amin_amax}, we have
\begin{equation*}
	\begin{aligned}
		\sum_{d^{k}\neq0}\|x^{k+1}-x^{k}\|^{2}
		=\sum_{d^{k}\neq0}\|t_{k}d^{k}\|^{2}
		=\sum_{d^{k}\neq0}\dfrac{\delta^{2}\psi^{2}(x^{k},d^{k})\|d^{k}\|^{2}}{\|d^{k}\|_{B^{k}}^{4}}
		\leq\dfrac{\delta^{2}}{a_{\min}^{2}}\sum_{d^{k}\neq0}\dfrac{\psi^{2}(x^{k},d^{k})}{\|d^{k}\|^{2}},
	\end{aligned}
\end{equation*}
which together with Lemma \ref{zou1} gives us 
$\sum_{d^{k}\neq0}\|x^{k+1}-x^{k}\|^{2}<\infty.$
This implies that $\|x^{k+1}-x^{k}\|^{2}\rightarrow0$, and thus $\lvert\beta_{k}^{\rm PRP}\rvert\rightarrow0$ and $\beta_{k}\rightarrow0$
as $k\rightarrow\infty$. By the definition of $d^{k}$, we have
\begin{equation}\label{prp_theorem0}
	\begin{aligned}
		\|d^{k}\|&=\|v(x^{k})+\beta_{k}d^{k-1}\|\\
		&\leq\|v(x^{k})\|+\beta_{k}\|d^{k-1}\|\\
		&\leq\|v(x^{k})\|+\beta_{k}\|v(x^{k-1})\|+\beta_{k}\beta_{k-1}\|d^{k-2}\|\\
		&\leq\cdots\\
		&\leq\|v(x^{k})\|+\beta_{k}\|v(x^{k-1})\|+\cdots+\beta_{k}\beta_{k-1}\cdots\beta_{1}\|v(x^{0})\|.
	\end{aligned}
\end{equation}
From Assumption {\normalfont\ref{A1}}, it follows that $\mathcal{L}$ is bounded, and thus $\{x^{k}\}\subset\mathcal{L}$ and $\{\|v(x^{k})\|\}$ is bounded using Proposition \ref{pa_sta_equ}\eqref{pa_sta_equ3}. Furthermore, we can obtain $\{\|d^{k}\|\}$ is bounded, and $\{\|JF(x^{k})\|_{2,\infty}\}$ is bounded above because the continuity argument for $JF$. 

Observe that
\begin{equation*}
	\lvert\psi(x^{k},d^{k-1})\rvert =\left\lvert\max_{i\in\langle m\rangle}\langle\nabla f_{i}(x^{k}),d^{k-1}\rangle\right\rvert\leq\max_{i\in\langle m\rangle}\|\nabla f_{i}(x^{k})\|\|d^{k-1}\|=\|JF(x^{k})\|_{2,\infty}\|d^{k-1}\|.
\end{equation*}
By the definition of $d^{k}$, the non-negativeness of $\beta_{k}$ and Lemma \ref{h_property}\eqref{h_property1}--\eqref{h_property2}, one has
\begin{equation*}
	\psi(x^{k},d^{k})\leq\psi(x^{k},v(x^{k}))+\beta_{k}\psi(x^{k},d^{k-1})\leq\psi(x^{k},v(x^{k}))+\beta_{k}\lvert\psi(x^{k},d^{k-1})\rvert,
\end{equation*}
which, together with the relation $ \psi(x^{k},v(x^{k}))\leq-\|v(x^{k})\|^{2}/2$, yields
\begin{equation}\label{prp_theorem1}
	\begin{aligned}
		-\psi(x^{k},d^{k})&\geq-\psi(x^{k},v(x^{k}))-\beta_{k}\lvert\psi(x^{k},d^{k-1})\rvert\\
		&\geq\dfrac{\|v(x^{k})\|^{2}}{2}-\beta_{k}\|JF(x^{k})\|_{2,\infty}\|d^{k-1}\|\\
		&\geq\dfrac{\|v(x^{k})\|^{2}}{4}\\
		&\geq\dfrac{\gamma^{2}}{4}>0
	\end{aligned}
\end{equation}
for large $k$, where the third inequality holds in view of $$\beta_{k}\|d^{k-1}\|\leq\dfrac{\|v(x^{k})\|^{2}}{4\|JF(x^{k})\|_{2,\infty}}$$ for large $k$, as $\beta_{k}\rightarrow0$ for large $k$ and $\{\|d^{k-1}\|\}$ has a bound. By \eqref{prp_theorem1} and the upper boundedness of $\|d^{k}\|$, we conclude that there exists $c_{1}>0$ such that
\begin{equation*}
	\begin{aligned}
		\sum_{d^{k}\neq0}\dfrac{\psi^{2}(x^{k},d^{k})}{\|d^{k}\|^{2}}\geq\sum_{d^{k}\neq0}\dfrac{\gamma^{2}}{4}\dfrac{1}{\|d^{k}\|^{2}}\geq\sum_{d^{k}\neq0}\dfrac{\gamma^{2}}{4c_{1}}=\infty.
	\end{aligned}
\end{equation*}
This contradicts Lemma \ref{zou1}. Therefore, the proof is complete.
\qed

\begin{theorem}\label{hs_theorem}
	Suppose that Assumptions {\normalfont\ref{A2}}--{\normalfont\ref{A4}} hold. Consider Algorithm {\normalfont\ref{alg1}} with $\beta_{k}=\max\{\beta_{k}^{\rm HS},0\}$. If $d^{k}$ is a descent direction, then $\lim\inf_{k\rightarrow\infty}\|v(x^{k})\|=0$.
\end{theorem}
\noindent{\bf Proof.}  
Assume by contradiction that there exists a positive constant $\gamma$ such that 
$$\|v(x^{k})\|\geq\gamma$$ for all $k\geq 0$. By \eqref{iteration}, the definition of $t_{k}$ and \eqref{amin_amax}, we have
\begin{equation*}
	\begin{aligned}
		\sum_{d^{k}\neq0}\|x^{k+1}-x^{k}\|^{2}
		=\sum_{d^{k}\neq0}\|t_{k}d^{k}\|^{2}
		=\sum_{d^{k}\neq0}\dfrac{\delta^{2}\psi^{2}(x^{k},d^{k})\|d^{k}\|^{2}}{\|d^{k}\|_{B^{k}}^{4}}
		\leq\dfrac{\delta^{2}}{a_{\min}^{2}}\sum_{d^{k}\neq0}\dfrac{\psi^{2}(x^{k},d^{k})}{\|d^{k}\|^{2}},
	\end{aligned}
\end{equation*}
which implies that
$\sum_{d^{k}\neq0}\|x^{k+1}-x^{k}\|^{2}<\infty$
because Lemma \ref{zou1}. Hence 
\begin{equation}\label{hs_theorem0_0}
	\|x^{k+1}-x^{k}\|^{2}\rightarrow0
\end{equation}
as $k\rightarrow\infty$. By the definition of $\beta_{k}^{\rm HS}$ and \eqref{lemma1_1}, we have
\begin{equation}\label{hs_theorem1_0}
	\begin{aligned}
		\beta_{k}^{\rm HS}
		=\dfrac{-\psi(x^{k},v(x^{k}))+\psi(x^{k-1},v(x^{k}))}{\psi(x^{k},d^{k-1})-\psi(x^{k-1},d^{k-1})}
		=\dfrac{-\psi(x^{k},v(x^{k}))+\psi(x^{k-1},v(x^{k}))}{(\rho_{k-1}-1)\psi(x^{k-1},d^{k-1})}.
	\end{aligned}
\end{equation}
Combining it with the definition of $d^{k}$ and Proposition \ref{h_property}\eqref{h_property1}--\eqref{h_property2}, we have
\begin{equation*}\label{hs_theorem1}
	\begin{aligned}
		\psi(x^{k},d^{k})&=\psi(x^{k},v(x^{k})+\beta_{k}d^{k-1})\\
		&\leq\psi(x^{k},v(x^{k}))+\beta_{k}\psi(x^{k},d^{k-1})\\
		&=\psi(x^{k},v(x^{k}))+\beta_{k}\rho_{k-1}\psi(x^{k-1},d^{k-1})\\
		&=\psi(x^{k},v(x^{k}))+\dfrac{\rho_{k-1}}{\rho_{k-1}-1}(-\psi(x^{k},v(x^{k}))+\psi(x^{k-1},v(x^{k}))).
	\end{aligned}
\end{equation*}
Dividing both sides of the above inequality by $\psi(x^{k},v(x^{k}))$, and observing that $\psi(x^{k},v(x^{k}))<0$ and \eqref{hs_theorem0_0}, we obtain
\begin{equation*}
	\begin{aligned}
		\dfrac{\psi(x^{k},d^{k})}{\psi(x^{k},v(x^{k}))}\geq 1+\dfrac{\rho_{k-1}}{\rho_{k-1}-1}\dfrac{-\psi(x^{k},v(x^{k}))+\psi(x^{k-1},v(x^{k}))}{\psi(x^{k},v(x^{k}))}\geq \dfrac{1}{2}
	\end{aligned}
\end{equation*}
for large $k$. Then, remembering that $\|v(x^{k})\|^{2}\leq -2\psi(x^{k},v(x^{k}))$, we have
\begin{equation*}\label{hs_theorem2}
	\begin{aligned}
		\dfrac{4\psi^{2}(x^{k},d^{k})}{\|v(x^{k})\|^{4}}\geq\dfrac{\psi^{2}(x^{k},d^{k})}{\psi^{2}(x^{k},v(x^{k}))}\geq \dfrac{1}{4}
	\end{aligned}
\end{equation*}
for large $k$, i.e., $\lvert\psi(x^{k},d^{k})\rvert\geq\|v(x^{k})\|^{2}/4$ for large $k$. This, together with \eqref{hs_theorem1_0} and Corollary \ref{coro1}\eqref{coro12}, yields
\begin{equation*}\label{hs_theorem1_1}
	\begin{aligned}
		\lvert\beta_{k}^{\rm HS}\rvert&=\left\lvert\dfrac{-\psi(x^{k},v(x^{k}))+\psi(x^{k-1},v(x^{k}))}{(1-\rho_{k-1})\psi(x^{k-1},d^{k-1})}\right\rvert\\
		&\leq\frac{4}{1-\rho_{k-1}}\dfrac{\lvert-\psi(x^{k},v(x^{k}))+\psi(x^{k-1},v(x^{k}))\rvert}{\|v(x^{k-1})\|^{2}}\\
		&\leq\frac{4a_{\max}}{\mu \delta}\dfrac{\lvert-\psi(x^{k},v(x^{k}))+\psi(x^{k-1},v(x^{k}))\rvert}{\|v(x^{k-1})\|^{2}}\\
		&\leq\frac{4a_{\max}}{\mu \delta\gamma^{2}}\lvert-\psi(x^{k},v(x^{k}))+\psi(x^{k-1},v(x^{k}))\rvert,
	\end{aligned}
\end{equation*}
which combined with \eqref{hs_theorem0_0} yields $\lvert\beta_{k}^{\rm HS}\rvert\rightarrow0$ and $\beta_{k}\rightarrow0$
as $k\rightarrow\infty$. From \eqref{prp_theorem0}, we have that $\{\|d^{k}\|\}$ is bounded. Then, the remaining proof is the same as that of Theorem \ref{prp_theorem}.
\qed

\section{Numerical experiments}\label{sec:5}

This section presents numerical experiments to demonstrate the performance of the proposed stepsize strategy. Our experiments are conducted using \texttt{MATLAB R2020b} on a computer with the following characteristics: 11th Gen Intel(R) Core(TM) i7-11390H processor (3.40 GHz) and 16 GB of RAM.

It is noteworthy that the proposed stepsize rule is governed by $\delta$ and $B^{k}$. In the following, we will begin by discussing how to update these two parameters during the iterations.

\subsection{The update techniques of $\delta$ and $B^{k}$}

In our numerical experiments, the case of $B^{k}=I$ for all $k\geq0$ shows unsatisfactory convergence behavior. Therefore, we consider another update way for $B^{k}$, i.e., $B^{0}=I$, and the rest iterations use the modified BFGS update formula:
\begin{equation}\label{bfgs}
	B^{k+1}=\left\{\begin{array}{lllll}
		B^{k}+\frac{y^{k}(y^{k})^{\top}}{(s^{k})^{\top}y^{k}}-\frac{B^{k}s^{k}(s^{k})^\top B^{k}}{(s^{k})^\top B^{k}s^{k}}, &\quad (s^{k})^{\top}y^{k}>0,\\
		B^{k}, & \quad{\rm otherwise},
	\end{array}\right.
\end{equation}
where
\begin{equation*}
	\begin{aligned}
		s^{k}&=x^{k+1}-x^{k},\\
		y^{k}&=\sum_{i=1}^{m}\lambda_{i}^{k}(\nabla f_{i}(x^{k+1})-\nabla f_{i}(x^{k})).
	\end{aligned}
\end{equation*}
It is noteworthy that $\lambda_{i}^{k}$ ($i\in\langle m\rangle$) is the optimal solution of \eqref{dual}, which slightly differs from the setting of \cite{ansary2015modified,lapucci2023limited,chen2022variable}. 

Since $\delta$ in \eqref{delta} depends on the unknown parameters $a_{\min}$ and $L$, 
determining its precise value in practice can be quite challenging. If the value of $\delta$ is artificially set too large in the test, the CG methods with our stepsize rule may occasionally produce an $x^{k+1}$ such that $F(x^{k+1})\npreceq F(x^{k}) $. In theory, as shown in \eqref{F_descent}, this issue can be avoided if $\delta$ is suffciently small. When the case of $F(x^{k+1})\npreceq F(x^{k}) $ is detected in practice, we reduce $\delta$ by half\footnote{The initial value of $\delta$ is set to 1.}. Thus, the method will generate a satisfactory iterate after finite steps and finally converge.  The technique for setting $\delta$ is similar to the approach used in \cite{chen2002global} for solving single objective optimization problems.

\subsection{Experiments settings}

In the following subsections, we present a detailed description of the experimental settings. 

\subsubsection{Algorithms}
In the classical multiobjective steepest descent (SD) method \cite{fliege2000steepest}, the stepsize was obtained by the Armijo line search. According to Remark \ref{re4.1}, the SD\_no\footnote{We use ``SD\_no" to denote the SD method \cite{fliege2000steepest} with the proposed stepsize formula \eqref{tk}.} method converges. Therefore, we first compare the SD\_no method with the SD method. Additionally, we also compare the CG methods using the stepsize \eqref{tk} with the CG methods using the Wolfe-type line search \cite{lucambio2018nonlinear}. For brevity, we use ``FR\_no" to denote the CG method with $\beta_{k}^{\rm FR}$ and the proposed stepsize formula. By doing this, we have five distinct algorithms, i.e., FR\_no, CD\_no, DY\_no, PRP\_no and HS\_no.

All the compared algorithms require solving the optimization problem \eqref{sub_pro}  to determine $v(x^{k})$. In the test, we consider the corresponding dual problem of \eqref{sub_pro} (see \cite{fliege2000steepest}):
\begin{equation}\label{dual}
	\begin{aligned}
		\min&\quad \dfrac{1}{2}\left\|\sum_{i=1}^{m}\lambda_{i}\nabla f_{i}(x^{k})\right\|^{2}\\
		{\rm s.t.}&\quad \lambda\in\varLambda^{m},
	\end{aligned}
\end{equation}
where $\varLambda^{m}=\{\lambda\in\mathbb{R}^{m}:\sum_{i=1}^{m}\lambda_{i}=1,\lambda_{i}\geq0,\forall i\in\langle m\rangle\}$ stands for the simplex. Then, $v(x^{k})$ can also be represented as
\begin{equation*}\label{dual_sol}
	v(x^{k})=-\sum_{i=1}^{m}\lambda_{i}^{k}\nabla f_{i}(x^{k}),
\end{equation*}
where $\lambda_{i}^{k}=(\lambda_{1}^{k},\lambda_{2}^{k},\ldots,\lambda_{m}^{k})\in \varLambda^{m}$ is optimal solution of $\eqref{dual}$. The standard MATLAB subroutine \texttt{quadprog} is adopted to solve the dual problem \eqref{dual}.

All runs are stopped whenever $\lvert\theta(x^{k})\rvert$ is less than or equal $10^{-6}$. Since Proposition \ref{pa_sta_equ} implies that $v(x) = 0$ if and only if $\theta(x) = 0$, this stopping criterion makes sense. The maximum number of allowed outer iterations is set to 10000. In the Armijo line search \cite{fliege2000steepest}, we set the parameter $\beta=10^{-4}$, and in the strong Wolfe line search \cite{lucambio2018nonlinear}, we set the parameters $\rho_{1}=10^{-4}$ and $\rho_{2}=10^{-1}$.
 
\subsubsection{Test problems}

We select 37 test problems, consisting of both convex and nonconvex multiobjective optimization problems, as shown in Table \ref{testpro}. The first column identifies the problem name. The second and third columns, labeled
as ``$n$'' and ``$m$'', respectively, inform the number of variables and the number of objective functions of the problem. The initial points are generated within a box defined by the lower and upper bounds, denoted as $x_{L}$ and $x_{U}$, respectively. The last column indicates the the corresponding reference. It is important to emphasize that the boxes presented in Table \ref{testpro} are exclusively used for specifying the initial points and are not considered by the algorithms during their execution. 

\begin{table}\footnotesize
	\centering
	\caption{Test problems.}
	\vskip-0.1in
	\setlength{\tabcolsep}{6mm}{
		\begin{tabular}{llllll}
			\toprule
			Problem & $n$     &$ m $    & $x_{L}$     & $x_{U}$     & Source \\\midrule
			BK1   & 2     & 2     & $(-5,-5)$ & $(10,10)$ & \cite{huband2006review} \\
			DD1a$^{1}$   & 5     & 2     & $(-1,\ldots,-1)$ & $(1,\ldots,1)$ & \cite{das1998normal} \\
			DD1b$^{1}$   & 5     & 2     & $(-10,\ldots,-10)$ & $(10,\ldots,10)$ & \cite{das1998normal} \\
			DD1c$^{1}$   & 5     & 2     & $(-20,\ldots,-20)$ & $(20,\ldots,20)$ & \cite{das1998normal} \\
			DGO1   & 1    & 2     & $-10$ & $13$ & \cite{huband2006review} \\
			Far1   & 2    & 2     & $(-1,-1)$ & $(1,1)$ & \cite{huband2006review} \\
			FDSa   & 10     & 3     & $(-2,\ldots,-2) $& $(2,\ldots,2)$ & \cite{fliege2009newton}\\	
			FDSb   & 200     & 3     & $(-2,\ldots,-2) $& $(2,\ldots,2)$ & \cite{fliege2009newton}\\	
			FDSc   & 500     & 3     & $(-2,\ldots,-2) $& $(2,\ldots,2)$ & \cite{fliege2009newton} \\
			FDSd   & 1000     & 3     & $(-2,\ldots,-2) $& $(2,\ldots,2)$ & \cite{fliege2009newton} \\
			FF1   & 2     & 2     & $(-1,-1) $& $(1,1)$ & \cite{huband2006review} \\
			Hil1   & 2     & 2     & $(0,0) $& $(1,1)$ & \cite{hillermeier2001generalized} \\
			IKK1   & 2     & 3    & $(-50,-50) $& $(50,50)$ & \cite{hillermeier2001generalized} \\
			IM1   & 2     & 2    & $(1,1) $& $(4,2)$ & \cite{hillermeier2001generalized} \\
			JOS1a 	& 50  & 2     & $(-100,\ldots,-100) $&$ (100,\ldots,100) $& \cite{jin2001dynamic}  \\
			JOS1b 	& 500  & 2     & $(-100,\ldots,-100) $&$ (100,\ldots,100) $& \cite{jin2001dynamic}  \\
			JOS1c 	& 1000  & 2     & $(-100,\ldots,-100) $&$ (100,\ldots,100) $& \cite{jin2001dynamic}  \\
			KW2 	& 2  & 2     & $(-3,-3) $&$ (3,3) $& \cite{kim2005adaptive}  \\
			Lov1  & 2     & 2     &$ (-10,-10) $&$ (10,10)$ & \cite{lovison2011singular} \\
			Lov3  & 2     & 2     &$ (-100,-100) $&$ (100,100)$ & \cite{lovison2011singular} \\
			Lov4  & 2     & 2     & $(-20,-20) $& $(20,20)$ &  \cite{lovison2011singular}\\
			Lov5  & 3     & 2     &$ (-2,-2,-2) $& $(2,2,2)$ & \cite{lovison2011singular} \\
			MGH16$^{2}$ & 4    & 5    & $(-25,-5,-5,-1) $& $(25,5,5,1)$ & \cite{more1981testing}  \\
			MGH26$^{2}$ & 4     & 4     & $(-1,-1,-1,-1)$ & $(1,1,1,1)$ &\cite{more1981testing}   \\
			MOP2  & 2     & 2     & $(-4,-4) $& $(4,4)$ &  \cite{huband2006review}\\
			MOP3  & 2     & 2     & $(-\pi,-\pi) $& $(\pi,\pi)$ &  \cite{huband2006review}\\
			MOP5  & 2     & 3     &$ (-30,-30)$ & $(30,30)$ & \cite{huband2006review} \\
			PNR   & 2     & 2     &$ (-2,-2) $& $(2,2)$ &  \cite{preuss2006pareto}\\
			MMR5a   & 50     & 2     &$ (-5,\ldots,-5) $& $(5,\ldots,5)$ &  \cite{miglierina2008box}\\
			MMR5b   & 200     & 2     &$ (-5,\ldots,-5) $& $(5,\ldots,5)$ &  \cite{miglierina2008box}\\
			MMR5c   & 500     & 2     &$ (-5,\ldots,-5) $& $(5,\ldots,5)$ &  \cite{miglierina2008box}\\
			SLCDT2 & 10    & 3     & $(-1,...,-1)$ & $(1,...,1)$ &\cite{schutze2008convergence}  \\
			SP1   & 2     & 2     &$ (-100,-100)$ & $(100,100)$ & \cite{huband2006review} \\
			SSFYY2    & 1     & 2     & $-100$ & $100$ & \cite{huband2006review} \\
			Toi9$^{2}$   & 4     & 4     & $(-1,-1,-1,-1) $& $(1,1,1,1)$ & \cite{toint1983test} \\
			Toi10$^{2}$   & 4     & 3     & $(-2,-2,-2,-2) $& $(2,2,2,2)$ & \cite{toint1983test} \\
			VU1 & 2     & 2     & $(-3,-3)$ & $(3,3)$ &\cite{huband2006review}    \\
			\bottomrule
	\end{tabular}}%
	\label{testpro}%
	\begin{tablenotes}\footnotesize
		\item[1]$^{1}$This is a modified version of DD1 (see \cite{mita2019nonmonotone} for more details).
		\item[2]$^{2}$This is an adaptation of a single objective optimization problem to the multiobjective setting (see \cite{mita2019nonmonotone} for more details).
	\end{tablenotes}
\end{table}%

Taking into account of numerical reasons, we use a scaled version of \eqref{mop}, which was proposed in \cite{goncalves2022study,goncalves2022globally}, i.e., 
\begin{equation}\label{mop1}
	\min_{x\in\mathbb{R}^{n}}\quad (r_{1} f_{1}(x),r_{2} f_{2}(x),...,r_{m}f_{m}(x))^{\top},
\end{equation}
where the scaling factors are computed as $r_{i}=1/\max\{1,\|\nabla f_{i}(x^{0})\|_{\infty}\}$, $i\in\langle m\rangle$, $x^{0}$ is the given initial point. As illustrated in \cite{goncalves2022study,hu2024alternative,goncalves2022globally}, \eqref{mop} is equivalent to \eqref{mop1}, since they have the same Pareto critical points.

\subsubsection{Metrics}\label{metric}
All problems are solved 100 times using initial points from a uniform random distribution inside a box specified in Table \ref{testpro}. The numerical results are reported by the median number of iterations (Iter), the median number of function evaluations (Feval), the median number of gradient evaluations (Geval) and the percentage of runs that have reached a critical point (\%).

To compare the ability of the algorithms to properly generate Pareto front approximations, we use the well-known \emph{Purity} and ($\Gamma$ and $\Delta$) \emph{Spread} metrics \cite{custodio2011direct}. Roughly speaking, the Purity metric measures the ability of the algorithm to find points
on the Pareto front of the problem, and the Spread metric measures the ability of the algorithm to obtain well-distributed points along the Pareto front.

In order to compare the performance of the tested algorithms more clearly, we use the performance profiles\footnote{The performance profiles in this paper are generated using the \texttt{MATLAB} code \texttt{perfprof.m}, which is freely available on the website \url{https://github.com/higham/matlab-guide-3ed/blob/master/perfprof.m}.} proposed  by Dolan and Mor\'{e} \cite{D_b2002}.

\subsection{Results analysis}

\begin{figure}[htbp]
	\centering   
	\includegraphics[width=\linewidth]{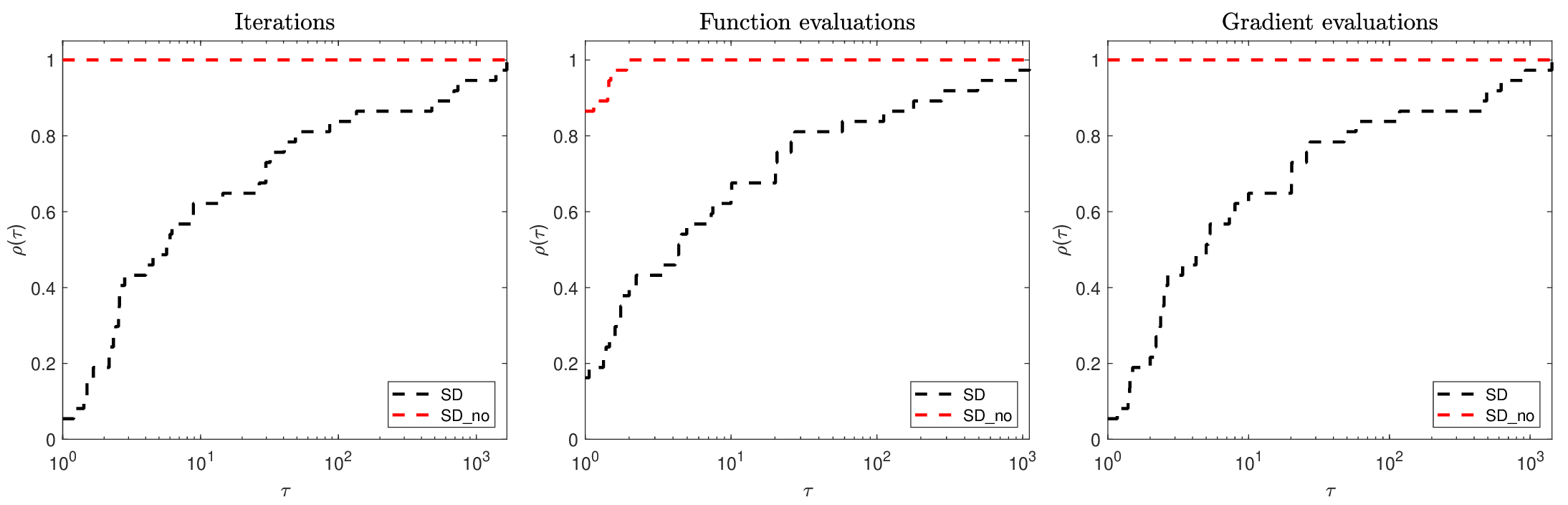}
	\caption{ Performance profiles using the number of iterations, the number of function evaluations and the number of gradient evaluations as the performance measurements for the SD and SD\_no methods.}
	\label{sd_it}
\end{figure}

\begin{figure}[htbp]
	\centering   
	\includegraphics[width=\linewidth]{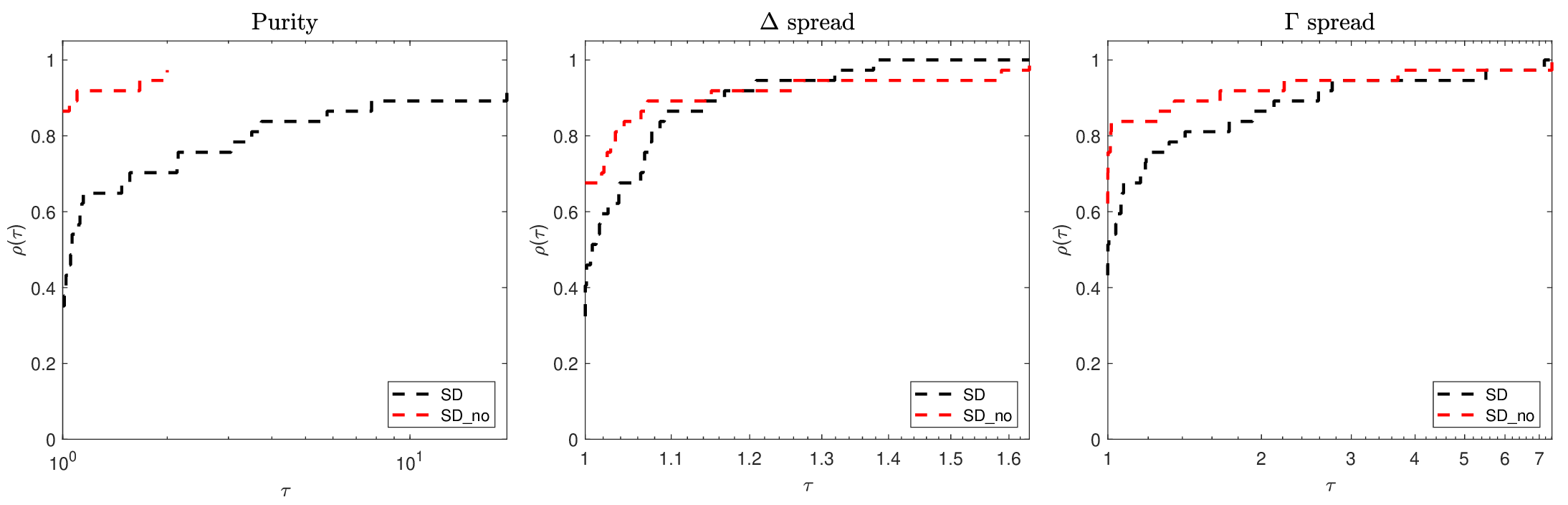}
	\caption{ Performance profiles using the Purity and ($\Delta$ and $\Gamma$) spread as the performance measurements for the SD and SD\_no methods.}
	\label{sd_sdno_per}
\end{figure}

In this section, we begin by comparing the performance of the SD method with that of the SD\_no method. Table \ref{sd_sdno_res} lists the performance results of the two algorithms on the benchmark instances. It can be seen that the SD method does not outperform the SD\_no method in terms of the number of iterations and gradient evaluations for all test problems. However, for function evaluations, the SD\_no method performs worse than the SD method on FDSa--FDSd and MOP2. For IM1 and MMR5c, the SD method struggles to converge within the maximum number of iterations, whereas the SD\_no method converges with fewer iterations. To offer a more intuitive representation of the numerical results, Fig. \ref{sd_it} presents the performance profiles for the number of iterations, function evaluations and gradient evaluations. As observed, the SD\_no method outperforms the SD method in finding stationary points on all test problems. Fig. \ref{sd_sdno_per} shows the Purity and ($\Delta$ and $\Gamma$) Spread performance profiles for the SD method and the SD\_no method. It can be observed that the SD\_no method outperforms the SD method in terms of the Purity and $\Gamma$ metrics, while the SD method shows better performance than the SD\_no method with respect to the $\Delta$ metric when $\tau>1.2$. Overall, the stepsize strategy we proposed performs better than the multiobjective Armijo line search strategy in the experimental setting of this paper.

\begin{table}[ht]\footnotesize
	\centering
	\caption{Numerical results of SD and SD\_no on the chosen set of test problems.}
	\begin{tabular}{lllllllll}
		\toprule
		& \multicolumn{4}{l}{SD}       & \multicolumn{4}{l}{SD\_no} \\
		\midrule
		Problem & Iter  & Feval & Geval & \%    & Iter  & Feval & Geval & \% \\
		\midrule
		BK1   & 29    & 60    & 60    & 100   & 2     & 6     & 6     & 100 \\
		DD1a  & 34    & 70    & 70    & 100   & 6     & 16    & 14    & 100 \\
		DD1b  & 30    & 62    & 62    & 100   & 5     & 15    & 12    & 100 \\
		DD1c  & 31    & 64    & 64    & 100   & 5     & 14    & 12    & 100 \\
		DGO1  & 2     & 6     & 6     & 100   & 2     & 6     & 6     & 100 \\
		Far1  & 18.5  & 67    & 39    & 100   & 13    & 63    & 28    & 100 \\
		FDSa  & 31    & 96    & 96    & 100   & 11    & 60    & 36    & 100 \\
		FDSb  & 106.5 & 322.5 & 322.5 & 100   & 44    & 483   & 135   & 100 \\
		FDSc  & 124.5 & 376.5 & 376.5 & 100   & 49    & 546   & 150   & 100 \\
		FDSd  & 135   & 408   & 408   & 100   & 52    & 582   & 159   & 100 \\
		FF1   & 21    & 44    & 44    & 100   & 9     & 30    & 20    & 100 \\
		Hil1  & 6     & 32    & 14    & 100   & 6     & 24    & 14    & 100 \\
		IKK1  & 59.5  & 181.5 & 181.5 & 100   & 2     & 9     & 9     & 100 \\
		IM1   & 10000 & 20002 & 20002 & 0     & 21    & 112   & 44    & 100 \\
		JOS1a & 173   & 348   & 348   & 100   & 2     & 6     & 6     & 100 \\
		JOS1b & 1472.5 & 2947  & 2947  & 100   & 2     & 6     & 6     & 100 \\
		JOS1c & 2775  & 5552  & 5552  & 100   & 2     & 6     & 6     & 100 \\
		KW2   & 10.5  & 41    & 23    & 100   & 7     & 23    & 16    & 100 \\
		Lov1  & 35.5  & 73    & 73    & 100   & 4     & 10    & 10    & 100 \\
		Lov3  & 80    & 162   & 162   & 100   & 3     & 8     & 8     & 100 \\
		Lov4  & 81    & 164   & 164   & 100   & 2     & 6     & 6     & 100 \\
		Lov5  & 71    & 150   & 144   & 100   & 8     & 20    & 18    & 100 \\
		MGH16 & 59.5  & 310   & 302.5 & 100   & 2     & 15    & 15    & 100 \\
		MGH26 & 11    & 48    & 48    & 100   & 5     & 24    & 24    & 100 \\
		MOP2  & 6     & 14    & 14    & 100   & 5     & 16    & 12    & 100 \\
		MOP3  & 15.5  & 35    & 33    & 100   & 6     & 20    & 14    & 100 \\
		MOP5  & 16    & 51    & 51    & 100   & 4     & 15    & 15    & 100 \\
		PNR   & 5     & 14    & 12    & 100   & 3     & 8     & 8     & 100 \\
		MMR5a & 2220.5 & 4443  & 4443  & 100   & 45.5  & 440   & 93    & 100 \\
		MMR5b & 6534  & 13070 & 13070 & 100   & 9.5   & 47    & 21    & 100 \\
		MMR5c & 10000 & 20002 & 20002 & 0     & 6     & 18    & 14    & 100 \\
		SLCDT2 & 18    & 57    & 57    & 100   & 7     & 25.5  & 24    & 100 \\
		SP1   & 945.5 & 1893  & 1893  & 100   & 7     & 17    & 16    & 100 \\
		SSFYY2 & 128   & 258   & 258   & 100   & 4     & 10    & 10    & 100 \\
		Toi9  & 9     & 50    & 40    & 100   & 6     & 36    & 28    & 100 \\
		Toi10 & 52    & 208.5 & 159   & 100   & 11.5  & 42    & 37.5  & 100 \\
		VU1   & 151   & 304   & 304   & 100   & 69.5  & 584   & 141   & 100 \\
		\bottomrule
	\end{tabular}%
	\label{sd_sdno_res}%
\end{table}%

\begin{figure}[htbp]
	\centering   
			\includegraphics[width=\linewidth]{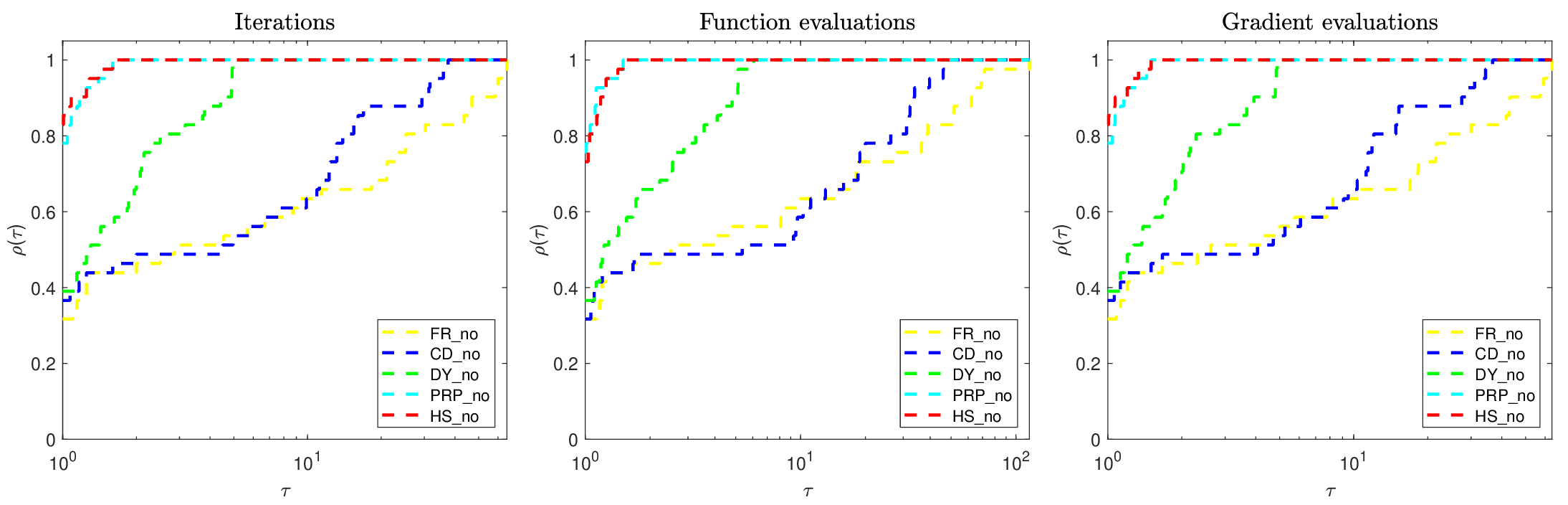}
	\caption{ Performance profiles using the number of iterations, the number of function evaluations and the number of gradient evaluations as the performance measurements for the FR\_no, CD\_no, DY\_no, PRP\_no and HS\_no methods.}
	\label{cgno_per}
\end{figure}

\begin{figure}[htbp]
	\centering   
		\includegraphics[width=\linewidth]{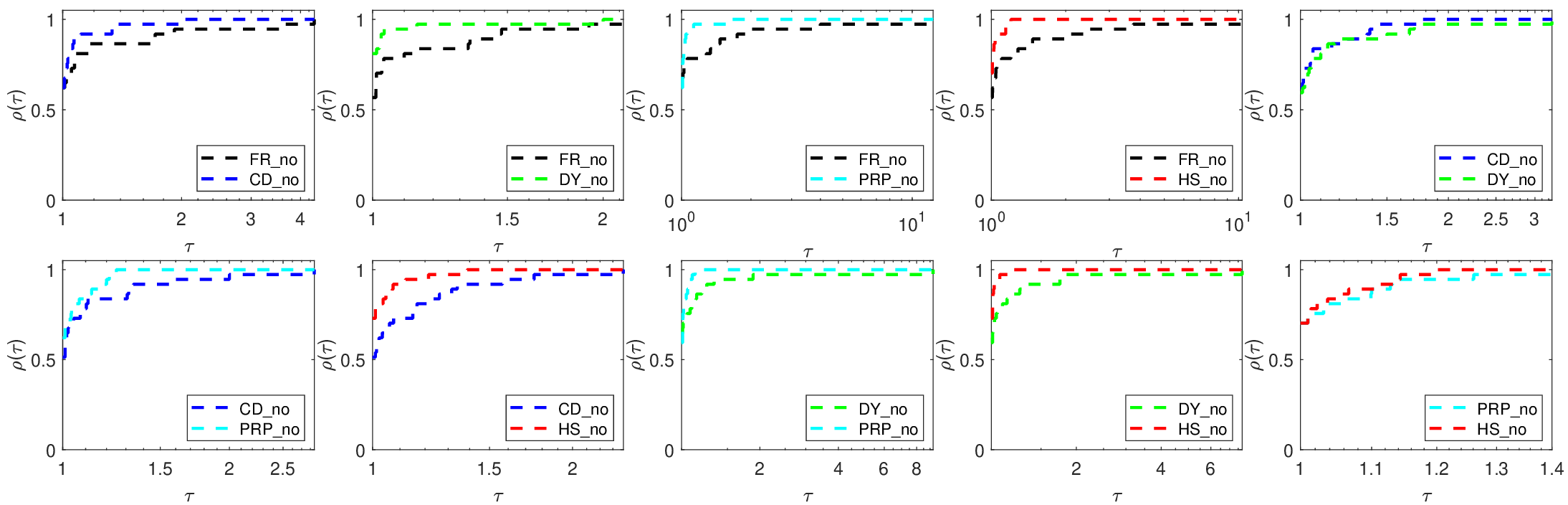}
	\caption{Purity performance profiles of the FR\_no, CD\_no, DY\_no, PRP\_no and HS\_no methods.}
	\label{cg_no_purity}
\end{figure}

\begin{figure}[htbp]
	\centering   
			\includegraphics[width=0.65\linewidth]{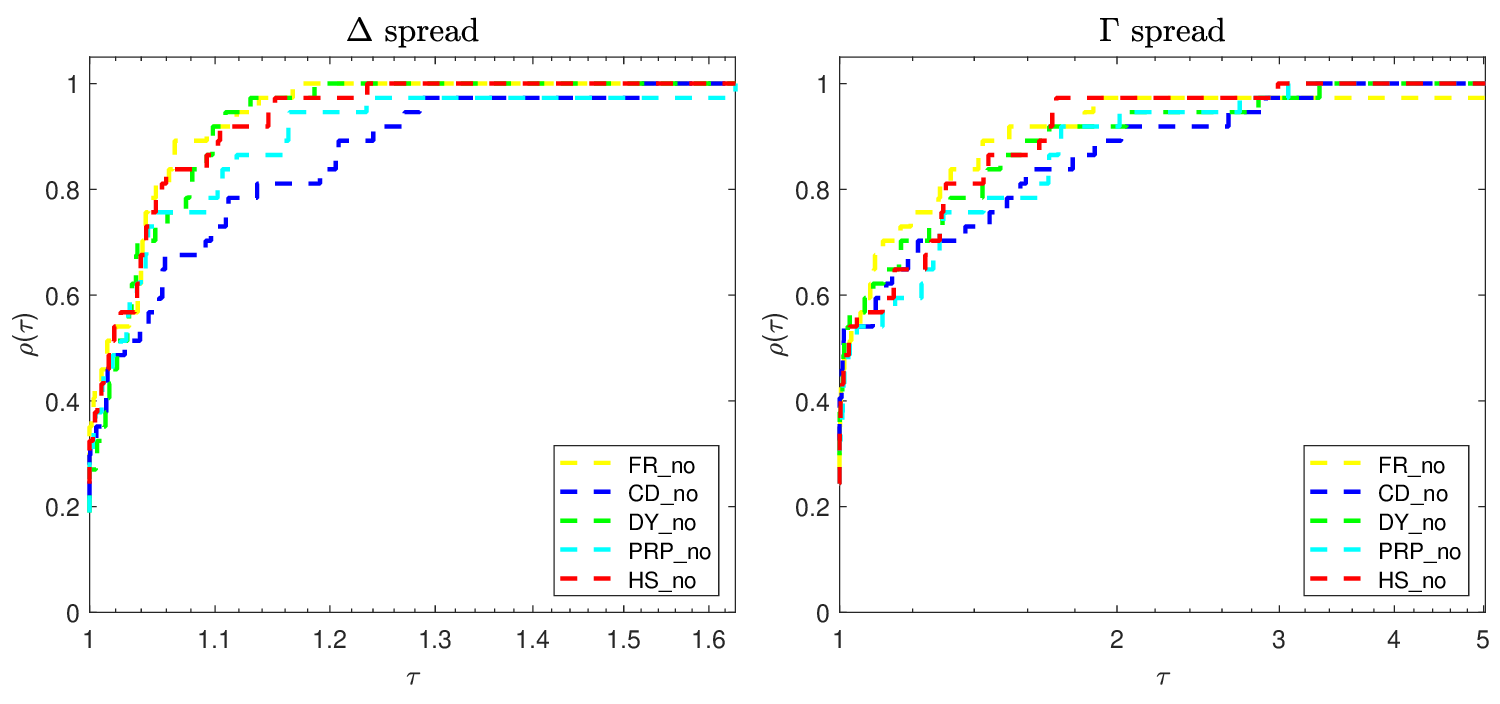}
	\caption{Spread performance profiles of the FR\_no, CD\_no, DY\_no, PRP\_no and HS\_no methods.}
	\label{cg_no_per}
\end{figure}

Next, we compare the FR\_no, CD\_no, DY\_no, PRP\_no and HS\_no methods to determine which parameter, $\beta_{k}$, is most effective under the proposed stepsize strategy. The numerical results of these five methods on the chosen set of test problems are listed in Table \ref{cg_no}. From this table, we can see that all five methods are able to reach critical points within the given maximum number of iterations. Fig. \ref{cgno_per} shows the results comparing the algorithms using as the performance measurement: the number of iterations, the number of function evaluations and the number of gradient evaluations. Considering the number of iterations,  the HS\_no method is the most effective (82.9\%), followed by the PRP\_no method (78.0\%), the DY\_no method (39.0\%), the CD\_no method (36.6\%) method and the FR\_no method (31.7\%). In terms of function and gradient evaluations, the HS\_no and the PRP\_no method clearly exhibit comparable performance, outperforming the rest three methods. Figs. \ref{cg_no_purity} and \ref{cg_no_per} show the Purity and ($\Delta$ and $\Gamma$) spread performance profiles for these algorithms. As recommended in \cite{custodio2011direct}, for the Purity metric, we compare the algorithms in pairs. For the Purity metric, it can be seen that HS\_no exhibits the best performance, followed closely by PRP\_no. For the $\Delta$ metric, there is no significant difference between FR\_no and HS\_no, both of which perform slightly better than the other methods. For the $\Gamma$ metric, no significant differences are observed among the considered algorithms.

\begin{sidewaystable}[thp]\footnotesize
	\centering
	\caption{Numerical results of the FR\_no, CD\_no, DY\_no, PRP\_no and HS\_no methods on the chosen set of test problems.}
	\resizebox{\linewidth}{!}{
	    \begin{tabular}{lllllllllllllllllllll}
		\toprule
		& \multicolumn{4}{l}{FR\_no}    & \multicolumn{4}{l}{CD\_no}    & \multicolumn{4}{l}{DY\_no}    & \multicolumn{4}{l}{PRP\_no}   & \multicolumn{4}{l}{HS\_no} \\
		\midrule
		Problem & Iter  & Feval & Geval & \%    & Iter  & Feval & Geval & \%    & Iter  & Feval & Geval & \%    & Iter  & Feval & Geval & \%    & Iter  & Feval & Geval & \% \\
		\midrule
		BK1   & 4     & 10    & 10    & 100   & 4     & 10    & 10    & 100   & 4     & 10    & 10    & 100   & 2     & 6     & 6     & 0     & 2     & 6     & 6     & 100 \\
		DD1a  & 151   & 448   & 304   & 100   & 83.5  & 497   & 169   & 100   & 26.5  & 97    & 55    & 100   & 6.5   & 16    & 15    & 0     & 6     & 16    & 14    & 100 \\
		DD1b  & 127   & 298   & 256   & 100   & 79    & 420   & 160   & 100   & 19    & 57    & 40    & 100   & 6     & 16    & 14    & 0     & 6     & 16    & 14    & 100 \\
		DD1c  & 126.5 & 269   & 255   & 100   & 79    & 516   & 160   & 100   & 15    & 41    & 32    & 100   & 6     & 16    & 14    & 0     & 6     & 16    & 14    & 100 \\
		DGO1  & 2     & 6     & 6     & 100   & 2     & 6     & 6     & 100   & 2     & 6     & 6     & 100   & 2     & 6     & 6     & 0     & 2     & 6     & 6     & 100 \\
		Far1  & 517   & 2068  & 1036  & 100   & 254   & 1102  & 510   & 100   & 42.5  & 196   & 87    & 100   & 24    & 117   & 50    & 0     & 20.5  & 114   & 43    & 100 \\
		FDSa  & 516   & 3091.5 & 1551  & 100   & 177   & 1125  & 534   & 100   & 41    & 246   & 126   & 100   & 11    & 60    & 36    & 0     & 11    & 60    & 36    & 100 \\
		FDSb  & 2574.5 & 29986.5 & 7726.5 & 100   & 1649.5 & 19188 & 4951.5 & 100   & 215   & 2475  & 648   & 100   & 44    & 483   & 135   & 0     & 44    & 483   & 135   & 100 \\
		FDSc  & 3042  & 35907 & 9129  & 100   & 1550.5 & 18483 & 4654.5 & 100   & 241   & 2797.5 & 726   & 100   & 49    & 546   & 150   & 0     & 49    & 546   & 150   & 100 \\
		FDSd  & 3390  & 40174.5 & 10173 & 100   & 1525.5 & 18646.5 & 4579.5 & 100   & 255   & 2973  & 768   & 100   & 52    & 582   & 159   & 0     & 52    & 582   & 159   & 100 \\
		FF1   & 114   & 310   & 230   & 100   & 122.5 & 640   & 247   & 100   & 21    & 71    & 44    & 100   & 10    & 32    & 22    & 0     & 10    & 33    & 22    & 100 \\
		Hil1  & 61    & 224   & 124   & 100   & 219   & 1241  & 440   & 100   & 13    & 49    & 28    & 100   & 7.5   & 29    & 17    & 0     & 7     & 27    & 16    & 100 \\
		IKK1  & 2     & 9     & 9     & 100   & 2     & 9     & 9     & 100   & 2     & 9     & 9     & 100   & 2     & 9     & 9     & 0     & 2     & 9     & 9     & 100 \\
		IM1   & 958.5 & 5764  & 1919  & 100   & 789.5 & 5336  & 1581  & 100   & 83.5  & 486   & 169   & 100   & 23    & 120   & 48    & 0     & 22    & 114   & 46    & 100 \\
		JOS1a & 2     & 6     & 6     & 100   & 2     & 6     & 6     & 100   & 2     & 6     & 6     & 100   & 2     & 6     & 6     & 0     & 2     & 6     & 6     & 100 \\
		JOS1b & 2     & 6     & 6     & 100   & 2     & 6     & 6     & 100   & 2     & 6     & 6     & 100   & 2     & 6     & 6     & 0     & 2     & 6     & 6     & 100 \\
		JOS1c & 2     & 6     & 6     & 100   & 2     & 6     & 6     & 100   & 2     & 6     & 6     & 100   & 2     & 6     & 6     & 0     & 2     & 6     & 6     & 100 \\
		KW2   & 45.5  & 135   & 93    & 100   & 99    & 528   & 200   & 100   & 13    & 36    & 28    & 100   & 10    & 28    & 22    & 0     & 10    & 33    & 22    & 100 \\
		Lov1  & 4     & 10    & 10    & 100   & 4     & 10    & 10    & 100   & 4     & 10    & 10    & 100   & 4     & 10    & 10    & 0     & 4     & 10    & 10    & 100 \\
		Lov3  & 3     & 8     & 8     & 100   & 3.5   & 9     & 9     & 100   & 3     & 8     & 8     & 100   & 3     & 8     & 8     & 0     & 3     & 8     & 8     & 100 \\
		Lov4  & 2     & 6     & 6     & 100   & 2     & 6     & 6     & 100   & 2     & 6     & 6     & 100   & 2     & 6     & 6     & 0     & 2     & 6     & 6     & 100 \\
		Lov5  & 7     & 16    & 16    & 100   & 7     & 16    & 16    & 100   & 7     & 16    & 16    & 100   & 7     & 18    & 16    & 0     & 7     & 16    & 16    & 100 \\
		MGH16 & 2     & 15    & 15    & 100   & 2     & 15    & 15    & 100   & 2     & 15    & 15    & 100   & 2     & 15    & 15    & 0     & 2     & 15    & 15    & 100 \\
		MGH26 & 5     & 24    & 24    & 100   & 4     & 26    & 20    & 100   & 5     & 24    & 24    & 100   & 5     & 24    & 24    & 0     & 5     & 28    & 24    & 100 \\
		MHHM2 & 2     & 9     & 9     & 100   & 2     & 9     & 9     & 100   & 2     & 9     & 9     & 100   & 3     & 168   & 12    & 0     & 2.5   & 10.5  & 10.5  & 100 \\
		MOP2  & 34    & 74    & 70    & 100   & 65.5  & 235   & 133   & 100   & 11    & 28    & 24    & 100   & 6     & 18    & 14    & 0     & 6.5   & 20    & 15    & 100 \\
		MOP3  & 20    & 49    & 42    & 100   & 31.5  & 118   & 65    & 100   & 10    & 26    & 22    & 100   & 7     & 22    & 16    & 0     & 7     & 23    & 16    & 100 \\
		MOP5  & 4     & 18    & 15    & 100   & 4     & 16.5  & 15    & 100   & 4     & 16.5  & 15    & 100   & 4     & 15    & 15    & 0     & 4     & 15    & 15    & 100 \\
		MOP7  & 9     & 30    & 30    & 100   & 7     & 24    & 24    & 100   & 8     & 27    & 27    & 100   & 8     & 27    & 27    & 0     & 7     & 24    & 24    & 100 \\
		PNR   & 5     & 12    & 12    & 100   & 5     & 12    & 12    & 100   & 4     & 12    & 10    & 100   & 4     & 10    & 10    & 0     & 4     & 10    & 10    & 100 \\
		MMR5a & 3302  & 43628 & 6606  & 100   & 837   & 11390 & 1676  & 100   & 234   & 3048  & 470   & 100   & 109   & 1198  & 220   & 0     & 116.5 & 1483  & 235   & 100 \\
		MMR5b & 664   & 7051  & 1330  & 100   & 72    & 680   & 146   & 100   & 20.5  & 175   & 43    & 100   & 14.5  & 61    & 31    & 0     & 14.5  & 69    & 31    & 100 \\
		MMR5c & 33.5  & 239   & 69    & 100   & 8     & 27    & 18    & 100   & 5     & 16    & 12    & 100   & 8     & 20    & 18    & 0     & 8     & 20    & 18    & 100 \\
		SLCDT2 & 73    & 255   & 222   & 100   & 135.5 & 1062  & 409.5 & 100   & 13    & 45    & 42    & 100   & 8     & 31.5  & 27    & 0     & 8     & 33    & 27    & 100 \\
		SP1   & 5.5   & 14    & 13    & 100   & 5     & 12    & 12    & 100   & 6     & 14    & 14    & 100   & 7     & 18    & 16    & 0     & 7     & 17    & 16    & 100 \\
		SSFYY2 & 4     & 10    & 10    & 100   & 4     & 10    & 10    & 100   & 4     & 10    & 10    & 100   & 4     & 10    & 10    & 0     & 4     & 10    & 10    & 100 \\
		Toi9  & 17.5  & 90    & 74    & 100   & 41    & 334   & 168   & 100   & 8     & 44    & 36    & 100   & 7     & 40    & 32    & 0     & 7     & 36    & 32    & 100 \\
		Toi10 & 8     & 27    & 27    & 100   & 7.5   & 25.5  & 25.5  & 100   & 7     & 24    & 24    & 100   & 10.5  & 36    & 34.5  & 0     & 9     & 36    & 30    & 100 \\
		VU1   & 4317.5 & 43007 & 8637  & 100   & 1112.5 & 11105 & 2227  & 100   & 355.5 & 2888  & 713   & 100   & 72    & 607   & 146   & 0     & 72.5  & 601   & 147   & 100 \\
		\bottomrule
	\end{tabular}%
	}
	\label{cg_no}%
\end{sidewaystable}%

\begin{figure}[htbp]
	\centering   
	\includegraphics[width=\linewidth]{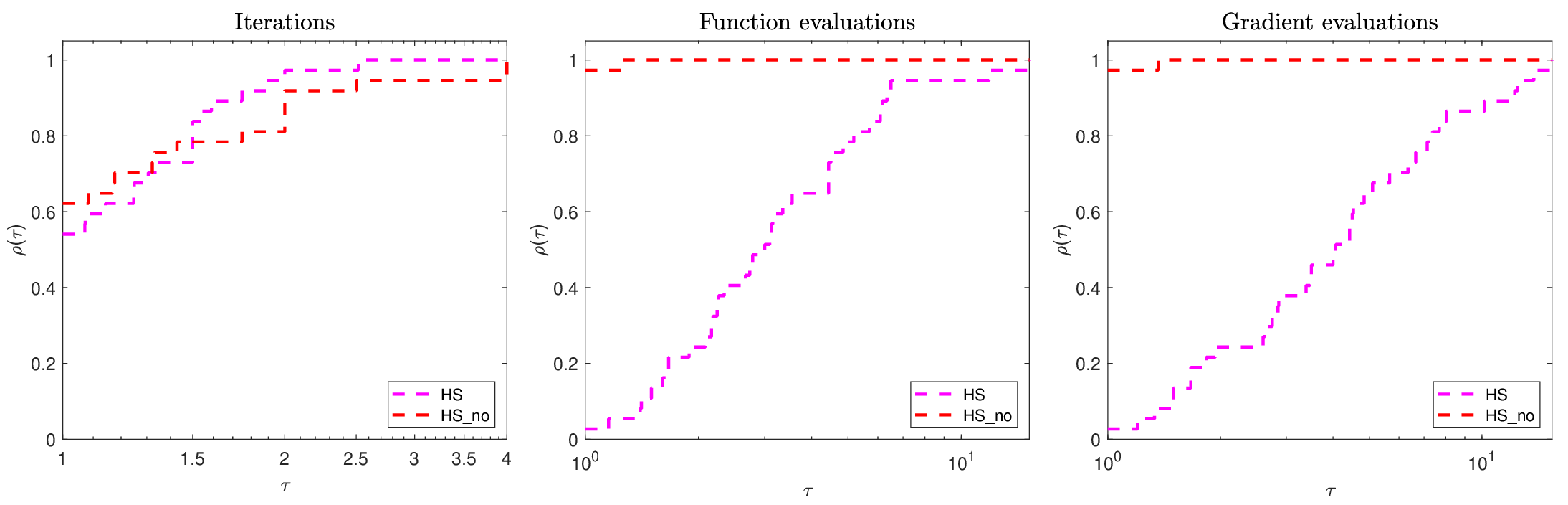}
	\caption{ Performance profiles using the number of iterations, the number of function evaluations and the number of gradient evaluations as the performance measurements for the HS and HS\_no methods.}
	\label{cg_per}
\end{figure}

\begin{figure}[htbp]
	\centering   
	\includegraphics[width=\linewidth]{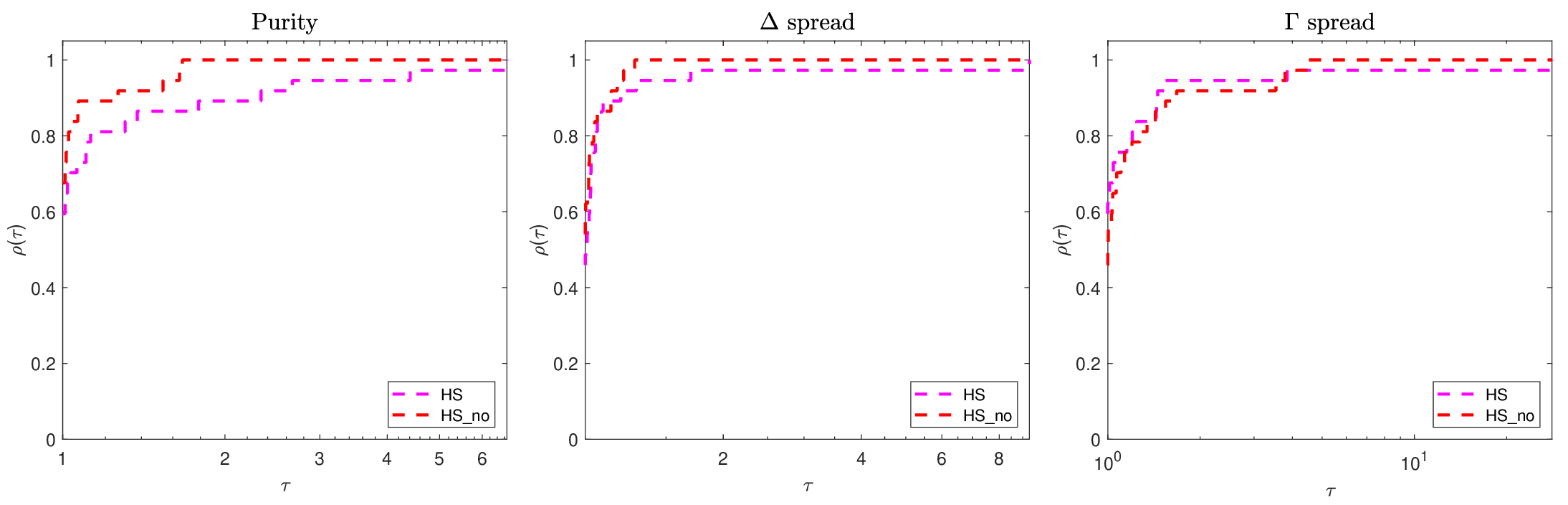}
	\caption{Purity and Spread performance profiles of the HS and HS\_no methods.}
	\label{hs_wolfe_no_metrics}
\end{figure}

It has been observed in \cite[pp. 2713]{lucambio2018nonlinear} that the HS and PRP methods, when applied with the strong Wolfe line search (i.e., \eqref{armijo} and \eqref{wolfeb}), are equivalent and exhibit superior performance compared to the FR, CD and DY methods under the same line search. From Fig. \ref{cgno_per}, it can be seen that the PRP\_no and HS\_no methods exhibit similar performance and show a certain advantage over the other three methods, namely FR\_no, CG\_no and DY\_no. Therefore, we next compare the proposed method, HS\_no, with the HS method \cite{lucambio2018nonlinear}. The numerical results obtained by the two algorithms for the selected instances are shown in Table \ref{hs_wolfe_no_res}.
When considering the iterations as an evaluation metric, the HS\_no method outperforms the HS method in the following 17 problem instances: BK1, DD1a--DD1c, FDSa--FDSd, Hil1, IM1, Lov4, MMR5a--MMR5c, Toi9, Toi10 and VU1. The HS\_no method is identical to the HS method in 6 problem instances: DGO1, Lov3, MGH16, PNR, SLCDT2 and SP1. In terms of function and gradient evaluations, the HS method outperforms the HS\_no method only on MOP5. The performance profiles obtained by the two algorithms are shown in Fig. \ref{cg_per}. From the first figure of Fig. \ref{cg_per}, it can be seen that when $\tau<1.5$, the HS\_no method show a significant advantage over the HS method. However, when $\tau\geq1.5$, the HS\_no method does not outperform the HS method. From the second and third figures of Fig. \ref{cg_per}, it is easy to see that the HS\_no method surpasses the HS method in function and gradient evaluations. The performance profiles of the Purity and Spread metrics are shown in Fig. \ref{hs_wolfe_no_metrics}.  As
can be seen, the HS\_no method is better than the HS method for the Purity and $\Delta$ metrics, whereas no signifcant diference is noticed for the $\Gamma$ metric.

\begin{table}[htbp]\footnotesize
	\centering
	\caption{Numerical results of the HS and HS\_no methods on the chosen set of test problems.}
	\begin{tabular}{lllllllll}
		\toprule
		& \multicolumn{4}{l}{HS}        & \multicolumn{4}{l}{HS\_no} \\
		\midrule
		Problem & Iter  & Feval & Geval & \%    & Iter  & Feval & Geval & \% \\
		\midrule
		BK1   & 4     & 39    & 34    & 100   & 2     & 6     & 6     & 100 \\
		DD1a  & 9     & 71    & 63    & 100   & 6     & 16    & 14    & 100 \\
		DD1b  & 9     & 71    & 62    & 100   & 6     & 16    & 14    & 100 \\
		DD1c  & 9     & 71    & 62    & 100   & 6     & 16    & 14    & 100 \\
		DGO1  & 2     & 13    & 11    & 100   & 2     & 6     & 6     & 100 \\
		Far1  & 15.5  & 131.5 & 113   & 100   & 20.5  & 114   & 43    & 100 \\
		FDSa  & 17    & 291   & 257   & 100   & 11    & 60    & 36    & 100 \\
		FDSb  & 55    & 1008  & 898   & 100   & 44    & 483   & 135   & 100 \\
		FDSc  & 64    & 1235  & 1107  & 100   & 49    & 546   & 150   & 100 \\
		FDSd  & 70    & 1362  & 1222  & 100   & 52    & 582   & 159   & 100 \\
		FF1   & 8.5   & 88    & 77    & 100   & 10    & 33    & 22    & 100 \\
		Hil1  & 7.5   & 83.5  & 72.5  & 100   & 7     & 27    & 16    & 100 \\
		IKK1  & 1     & 17    & 15    & 100   & 2     & 9     & 9     & 100 \\
		IM1   & 35    & 703   & 633   & 100   & 22    & 114   & 46    & 100 \\
		JOS1a & 1     & 9     & 8     & 100   & 2     & 6     & 6     & 100 \\
		JOS1b & 1     & 10    & 9     & 100   & 2     & 6     & 6     & 100 \\
		JOS1c & 1     & 10    & 9     & 100   & 2     & 6     & 6     & 100 \\
		KW2   & 7     & 74    & 63    & 100   & 10    & 33    & 22    & 100 \\
		Lov1  & 3     & 30    & 26    & 100   & 4     & 10    & 10    & 100 \\
		Lov3  & 3     & 25    & 22    & 100   & 3     & 8     & 8     & 100 \\
		Lov4  & 3     & 27    & 24    & 100   & 2     & 6     & 6     & 100 \\
		Lov5  & 4     & 35    & 31    & 100   & 7     & 16    & 16    & 100 \\
		MGH16 & 2     & 85.5  & 76.5  & 100   & 2     & 15    & 15    & 100 \\
		MGH26 & 2     & 45    & 40    & 100   & 5     & 28    & 24    & 100 \\
		MOP2  & 6     & 71    & 61    & 100   & 6.5   & 20    & 15    & 100 \\
		MOP3  & 6     & 63    & 55.5  & 100   & 7     & 23    & 16    & 100 \\
		MOP5  & 1     & 12    & 11    & 100   & 4     & 15    & 15    & 100 \\
		PNR   & 4     & 33.5  & 28.5  & 100   & 4     & 10    & 10    & 100 \\
		MMR5a & 125   & 2091  & 1889.5 & 100   & 116.5 & 1483  & 235   & 100 \\
		MMR5b & 36.5  & 420   & 379.5 & 100   & 14.5  & 69    & 31    & 100 \\
		MMR5c & 14    & 303.5 & 277   & 100   & 8     & 20    & 18    & 100 \\
		SLCDT2 & 8     & 104   & 91.5  & 100   & 8     & 33    & 27    & 100 \\
		SP1   & 7     & 88    & 77.5  & 100   & 7     & 17    & 16    & 100 \\
		SSFYY2 & 1     & 14    & 12    & 100   & 4     & 10    & 10    & 100 \\
		Toi9  & 8     & 228.5 & 203   & 100   & 7     & 36    & 32    & 100 \\
		Toi10 & 17    & 426.5 & 374   & 100   & 9     & 36    & 30    & 100 \\
		VU1   & 90.5  & 1675.5 & 1494.5 & 100   & 72.5  & 601   & 147   & 100 \\
		\bottomrule
	\end{tabular}%
	\label{hs_wolfe_no_res}%
\end{table}%

In summary, our numerical results on the selected benchmark problems indicate that, without consuming the number of objective function evaluations, the proposed stepsize strategy outperforms the multiobjective Armijo line search \cite{fliege2000steepest} and is competitive with multiobjective Wolfe-type line searches \cite{lucambio2018nonlinear}.

\section{Conclusion and remarks}\label{sec:6}

In this paper, we proposed a new stepsize strategy defined by the formula \eqref{tk} in multiobjective optimization, and successfully applied it to the CG methods proposed in \cite{lucambio2018nonlinear}. The Zoutendijk-type condition for multiobjective optimization are establised when using the proposed fixed stepsize formula. Our results reveal an interesting property of the CG methods for some choices of $\beta_{k}$, i.e., the global convergence can be guaranteed under mild assumptions by using the proposed stepsize rule
rather than the Wolfe-type conditions. Numerical experiments on a set of test problems indicate that the proposed stepsize scheme outperforms the multiobjective Armijo line search and is competitive with the Wolfe-type conditions.

As part of our future research agenda, we would like to further construct alternative way for updating the matrix $B^{k}$ and investigate the adaptive updating way for the parameter $\delta$. Secondly, we aim to apply the proposed stepsize rule to other variants of multiobjective conjugate gradient methods \cite{goncalves2020extension,goncalves2022study,yahaya2023hybrid}. Additionally, we hope that stepsize technique can also be useful for other multiobjective descent algorithms, such as multiobjective quasi-Newton algorithm \cite{prudente2022quasi,ansary2015modified,lapucci2023limited}. In scalar optimization case, Dai \cite{dai2011convergence} showed that  the global convergence properties for some conjugate gradient methods with constant stepsize, which opens up an intriguing new avenue for theoretical research in multiobjective optimization.

%
%

\section*{References}

\end{document}